\documentclass[11pt,reqno]{amsart}

\usepackage{amsmath}\allowdisplaybreaks  
\usepackage{amssymb}        
\usepackage{amsfonts}       
\usepackage{amsthm}         
\usepackage{mathrsfs}       
\usepackage{cases}          
\usepackage{siunitx}        

\usepackage{array}
\usepackage{booktabs}       
\usepackage{longtable}      
\usepackage{multirow}
\usepackage{makecell}       
\usepackage{tabularx}       
\usepackage{nicematrix}     

\usepackage{tikz}\usetikzlibrary{tikzmark, arrows.meta}   
\usetikzlibrary{positioning}
\usepackage{float} 
\usepackage{xcolor}           
\usepackage{graphicx}         
\usepackage{subcaption}

\usepackage[linesnumbered, ruled, vlined]{algorithm2e}  

\usepackage[colorlinks=true, linkcolor=blue, citecolor=blue, urlcolor=blue]{hyperref} 
\usepackage[nameinlink]{cleveref}	

\usepackage[left=2cm, right=2cm, top=2cm, bottom=2cm]{geometry}

\usepackage{bm} 
\usepackage[normalem]{ulem}           
\usepackage{lineno}         

\makeatletter
\newcommand\figcaption{\def\@captype{figure}\caption} 
\newcommand\tabcaption{\def\@captype{table}\caption}
\newcommand{\appref}[1]{\hyperref[#1]{Appendix~\ref*{#1}}}
\makeatother

\newtheorem{theorem}{Theorem}[section]          
\newtheorem{lemma}[theorem]{Lemma}              
\newtheorem{proposition}[theorem]{Proposition}  
\newtheorem{example}[theorem]{Example}          

\begin{document}
	\title[SCNN representation of FE function]{Sparsely connected neural network representation of Lagrange finite element function}
	
	\author[J.~Hao, Y.~Huang, N.~Yi]{Jiaxiong Hao$^{\dag}$, Yunqing Huang$^\S$, Nianyu Yi$^\dag$}
	\address{$\dag$ Hunan Key Laboratory for Computation and Simulation in Science and Engineering, School of Mathematics and Computational Science, Xiangtan University, Xiangtan 411105, Hunan, P.R.China} \email{moodbear@qq.com (J. Hao);\  yinianyu@xtu.edu.cn (N. Yi)}
	\address{$\S$ National Center for Applied Mathematics in Hunan, Key Laboratory of Intelligent Computing \& Information Processing of Ministry of Education, Xiangtan University, Xiangtan 411105, Hunan, P.R.China} \email{huangyq@xtu.edu.cn}
	
	\begin{abstract}
    We construct a mesh-induced sparsely connected neural network framework that exactly reproduces arbitrary-order Lagrange finite element spaces over simplicial meshes. Unlike conventional black-box neural surrogates, the proposed network architecture is fully dictated by finite element discretization rules: local computations stem from simplex geometry and barycentric coordinate transformations, while global consistency is enforced through shared degrees of freedom. 
    For linear Lagrange elements, local basis functions are directly implemented via affine barycentric layers, and high-order polynomial bases are explicitly decomposed into barycentric product compositions realized by specially designed $\mathrm{ReLU}^p$ modules. Equipped with element indicator branches and multiplication units, these modular local components are globally assembled into a sparsely connected neural network whose function space coincides exactly with the target finite element space, thereby inheriting the complete classical finite element approximation theory. By prescribing customized backward differentiation rules for piecewise activations, function values and their spatial gradients can be simultaneously extracted via automatic differentiation within a unified computational graph, eliminating the separate gradient calculation subroutines required in standard finite element implementations. Numerical experiments verify the accuracy of the neural network representation of Lagrange finite elements. Furthermore, by virtue of the intrinsic mesh-free nature of this neural network representation, finite element functions can be interpolated between non-matching meshes, and the proposed scheme can be applied to adaptive finite element methods for solving parabolic partial differential equations. An open-source code implementation of the proposed architecture is made publicly available \footnotemark.
	\end{abstract}
	
	\footnotetext{GitHub: \url{https://github.com/FEMmaster/Sparse-NN-Representation-of-FE-Space}.}
	
	
	
	\maketitle
\section{Introduction}
The expressive power of neural-network function classes is commonly studied in the context of fully connected architectures and universal approximation theory \cite{cybenko1989, hornik1989}. Although such results establish the theoretical possibility of approximating broad classes of functions, fully connected architectures are intrinsically globally coupled: each hidden unit acts over the entire input domain, and variations in individual parameters can influence the network response throughout the whole domain. 
Sparse connectivity provides a structured way to reduce this global coupling. Radial basis function (RBF) and kernel-based networks \cite{buhmann2003, wendland2004} exemplify this strategy, restricting interactions to local neighborhoods or specific centers. By imposing structural constraints on neural-network architectures \cite{Mocanu2018, Robinett2018, Stephenson2020}, such sparse designs reduce computational cost and improve scalability. Yet a critical limitation persists across all these paradigms: universal approximation theorems guarantee only the existence of an approximating network, without offering a constructive algorithmic pathway for deriving the network scale required to achieve a prescribed approximation accuracy. This gap motivates architectures with structured connectivity, in which the organization of parameters follows explicit generative rules that induce a function space equipped with a structured basis.
This philosophy is realized, in a rigorous mathematical form, by the classical discretization methods of scientific computing. The finite element method, for instance, provides such a rule by constructing approximation spaces from mesh resolution, polynomial degree, and degree-of-freedom sharing. Local support, sparse algebraic systems, and structured coupling between neighboring elements arise directly from this construction \cite{brenner2008, ern2004}. This structure suggests a direct route for transferring finite element approximation principles into neural-network design. Indeed, a finite element space may be viewed as a sparse computational graph in which local operations are determined by polynomial basis evaluation on individual elements, while the coupling between local components is governed by shared degrees of freedom. This perspective provides a way to transfer the locality, sparsity, and degree-of-freedom sharing of classical discretization methods into neural-network architectures, replacing heuristic connectivity choices with a discretization-induced structure.

Several recent works have explored the connection between finite element approximations and neural-network architectures from different perspectives. A rigorous connection between deep ReLU networks and finite element spaces was first established in \cite{He2020}, where it was proved that continuous piecewise linear finite element functions can be represented exactly by ReLU networks, thereby placing the relationship between ReLU architectures and linear finite elements on a precise theoretical footing. The representation, however, is built from global finite element basis functions, whose ReLU realization depends on the local mesh connectivity around each vertex. Consequently, the resulting architecture is patch-dependent rather than element-uniform, and its use on large-scale meshes is not accompanied by a concrete implementation procedure, numerical validation, or practically deployable construction strategy. In its original form, the construction also does not extend directly to nonconvex domains. An element-based viewpoint was adopted in \cite{Longo2023}, where a neural-network finite element construction with a unified architectural form was proposed. This approach alleviates, to some extent, the structural fragmentation present in earlier vertex-patch-based representations. However, the formulation does not retain basis functions and degrees of freedom as the primary organizing objects and therefore departs from the standard finite element description. Moreover, its network construction is highly elaborate, involving intricate compatibility constraints and elementwise network design, while concrete implementation details, numerical validation, and a readily implementable realization procedure are not provided.
A different line of work focuses on the neural-network realization of finite element shape functions. Zhang et al.~\cite{Zhang2020} proposed a structured representation that decomposes the algebraic expressions of shape functions into elementary operations, including piecewise linear mappings, multiplications, and inversions. 
This construction interprets interpolation functions as explicit computational graphs and provides a useful perspective on the algebraic realizability of finite element bases. However, the higher-dimensional construction is essentially based on elementwise tensor-product factorizations, which makes the formulation well suited to quadrilateral and hexahedral elements but not readily extensible to simplex elements. 
\v{S}kardov\'a et al.~\cite{kardov2025} further developed the HiDeNN framework and presented two implementations of finite element neural-network interpolation. The first implementation introduces an assembly layer that maps element-level auxiliary functions to global shape functions through a sparse affine transformation. Although this provides an explicit network representation of assembled shape functions, the assembly layer must be tailored to the element type and interpolation order. Moreover, the explicit construction of global shape functions is not required for the pointwise evaluation of a finite element interpolation, since the evaluation can be carried out on the element containing the query point using the corresponding local degrees of freedom. The second implementation adopts this local-evaluation perspective through a reference-element formulation, in which the query point is supplied together with the index of its containing element. Consequently, point localization is handled externally by a conventional mesh search algorithm, while the network only evaluates the local shape functions associated with the selected element.
 More recently, several works have explored the interface between finite element discretizations and neural-network representations from related but distinct perspectives. Jin~\cite{Jin2026} studies the relation between two-hidden-layer ReLU networks and finite elements, clarifying the representational overlap between shallow ReLU architectures and finite element spaces. Li et al.~\cite{Li2026} propose unfitted finite element interpolated neural networks, in which finite element interpolation is incorporated into neural-network approximations on unfitted meshes. Wang et al.~\cite{Wang2025NNEM} develop a neural network element method for PDEs, emphasizing elementwise neural-network structures within a numerical discretization framework. For broader developments on neural-network discretizations, finite-element-inspired architectures, and approximation properties of neural networks in finite element settings, we refer to \cite{DabySeesaram2025, He2023, Liu2023, Opschoor2020, Opschoor2024b, Saha2021}.
 
 This work develops a mesh-induced sparsely connected neural network framework to exactly represent Lagrange finite element spaces on simplicial meshes. Local network blocks for linear and high-order finite elements are explicitly built using barycentric coordinate transformations and ReLU-based product modules, and global finite element expressions are assembled via element indicator networks. For any Lagrange finite element function \(u_h\in V_h^p\) defined on a simplicial mesh \(\mathcal{T}_h\), we construct a neural network function \(u_{\theta}\) to globally represent \(u_h\). The architecture of this neural network is
\[
N_{\mathrm{nn}}^{(p)}
=
\left[
d,\,
\bigl(c(d+1),\sigma_1\bigr),\,
\bigl(c\,2^p n_{\mathrm{loc}},\sigma_p\bigr),\,
c\,n_{\mathrm{loc}},\,
(2c,\sigma_c),\,
\bigl(4c,\sigma_2\bigr),\,
1
\right],
\]
where \(c\) and \(n_{\mathrm{loc}}\) denote the number of elements and local degrees of freedom, respectively, and \((n,\sigma)\) denotes a hidden layer with \(n\) neurons and activation function \(\sigma\). By specifying customized backward differentiation rules for the nonlinear activations, function values and their gradients can be evaluated uniformly through automatic differentiation on a single computational graph. The resulting network space is mathematically identical to the target finite element space and therefore inherits the standard approximation theory and convergence properties of finite element methods.

Overall, our contributions are summarized as follows:
\begin{enumerate}
    \item For every function \(u_h\) in the Lagrange finite element space associated with the simplicial mesh \(\mathcal{T}_h\), we establish a mesh-induced sparsely connected neural network architecture that represents such finite element functions exactly.
    
    \item For a given simplicial finite element mesh, the geometric information of the mesh elements is encoded into the weights and biases between the input layer and the first hidden layer, thereby establishing a mapping between the associated Lagrange finite element space and the mesh-induced sparsely connected neural network. Consequently, the approximation theory of the sparsely connected neural network follows naturally from that of Lagrange finite element spaces.
    
    \item By designing suitable activation functions for the hidden-layer neurons, the derivatives of the represented finite element functions can be computed via automatic differentiation.
    
    \item Owing to the intrinsic mesh-free nature of the proposed representation framework, finite element functions can be exactly interpolated between non-matching meshes, providing a new implementation strategy for adaptive finite element algorithms for parabolic partial differential equations.
\end{enumerate}

The remainder of this paper is organized as follows. In Section \ref{secFENN}, we systematically construct the proposed sparsely connected neural network representation framework for Lagrange finite element spaces defined on simplicial meshes. Section \ref{secApp} numerically verifies the accuracy of the neural-network-based representation of finite element functions and further applies this technique to the adaptive finite element algorithm for parabolic equations. Section \ref{secCon} concludes the paper and outlines promising directions for future work, including the extension of the proposed neural representation paradigm to other types of finite element spaces and the development of efficient numerical algorithms based on the established network structure.

\section{Mesh-induced sparsely connected neural network}\label{secFENN}
    Starting from recall the definition of standard simplicial Lagrange finite elements and barycentric coordinate formulations, we first establish the network realization for linear finite elements by designing barycentric mapping layers, element indicator branches and multiplication modules, which enables exact local-to-global assembly of linear finite element neural representations. 
    We further extend the proposed architecture to arbitrary high-order Lagrange finite elements, where high-degree polynomial basis functions are decomposed into products of barycentric coordinates and implemented via dedicated $ReLU^p$ product modules; concrete examples including one-dimensional quadratic elements and two-dimensional cubic basis functions are provided to illustrate the detailed local network construction. 
    After finishing network structural design, we investigate the differentiability of the resulting sparsely connected neural network, specify modified backward derivative rules for custom activation functions, and realize automatic gradient evaluation of finite element fields relying solely on backpropagation without extra finite element gradient subroutines. 
    
\subsection{Lagrange finite element spaces on simplicial mesh}
    Let \(\mathcal T_h\) be a shape regular and simplicial mesh of \(\Omega\subset\mathbb{R}^d\), and $\mathcal{N}_h=\{x_i\}_{i=1}^N$ denotes its vertices. 
    
    For each element $K\in \mathcal{T}_h$, let $P_p(K)$ be the space of polynomials with total degree at most \(p\) on element \(K\). 
    The associated $p$th-order Lagrange finite element space is defined by
    \[
    V_h^p = \{v_h\in C(\Omega): v_h|_K\in P_p(K), \forall\, K\in\mathcal T_h\}.
    \]
    Let $\{z_i\}_{i=1}^{N_{\mathrm{dof}}}$ denote the global Lagrange nodes, where $N_{\mathrm{dof}}$ is the total number of degrees of freedom of the finite element space $V_h^p$. 
    The corresponding global Lagrange basis functions \(\{\Phi_i\}_{i=1}^{N_{\mathrm{dof}}}\subset V_h^p\) are characterized by the nodal interpolation property
    \[
    \Phi_i(z_j)=\delta_{ij},
    \qquad i,j=1,\dots,N_{\mathrm{dof}}.
    \]
    For each finite element approximation $u_h \in V_h^p$, 
    \begin{equation}\label{eq:global}
        u_h(x) = \sum_{i=1}^{N_{\mathrm{dof}}} u_i \Phi_i(x), \quad x \in \Omega,
    \end{equation}
    where $u_i \in \mathbb{R}$ is the value of $u_h(z_i)$ at the node associated with the $i$th Lagrange basis function.  
    The expression given in \eqref{eq:global} adopts a global formulation, while each basis function \(\Phi_i\) is assembled from piecewise polynomials defined over individual elements within its support. 
    From a computational perspective, the local nature of finite-element discretization is manifested via element-wise assembly of the variational formulation. 
    To convert the global formulation into element-level computation, an index mapping between local degrees of freedom on each element \(K\in\mathcal T_h\) and the global counterparts in \eqref{eq:global} must be established. 

    For each element \(K\in\mathcal T_h\), we define the local-to-global index map
    \[
    g_K:\{i\}_{i=1}^{n_{\mathrm{loc}}}\to\{g_{K}(i)\}_{i=1}^{n_{\mathrm{loc}}}\subset \{1,2,\cdots,N_{dof}\},
    \]
    where \(n_{\mathrm{loc}}=\dim P_p(K)=\binom{p+d}{d}\) is the number of local degrees of freedom on \(K\). 
    The local Lagrange basis $\{\phi_i\}_{i=1}^{n_{\mathrm{loc}}}\subset P_p(K)$ satisfies
    \[
    \phi_i = \Phi_{g_K(i)}|_K, \qquad i=1,\dots,n_{\mathrm{loc}}.
    \]  
    The finite element function \(u_h\) on \(K\) is then expressed through the local-to-global map \(g_K\) as
    \begin{equation}\label{eq:local}
        \left.u_h\right|_K(x) = \sum_{i=1}^{n_{\mathrm{loc}}} u_{g_K(i)}\phi_i(x), \qquad x\in K.
    \end{equation}

From the elementwise representation \eqref{eq:local}, the finite element solution \(u_h\) is obtained by superposing the local contributions of all elements. This assembly constitutes the standard procedure for constructing the global discrete algebraic system in conventional finite element frameworks. The same local-to-global structure is used to define finite element sparse neural networks, in which locality is inherited from the elementwise decomposition rather than imposed by an external sparsification rule.

    Let
    \[
    K=\operatorname{conv}\{a_1,\dots,a_{d+1}\}\subset\mathbb R^d
    \]
    be a nondegenerate $d$-simplex with affinely independent vertices $a_1,\dots,a_{d+1}$ with coordinates $\{x^{a_i}\}_{i=1}^{d+1}$. 
    Denote by
    \[
    K_i(x):=\operatorname{conv}\{a_1,\dots,a_{i-1}, a_x,a_{i+1},\dots,a_{d+1}\}
    \]
    the sub-simplex obtained from $K$ by replacing $a_i$ with $a_x$, which denotes the points with coordinate $x$. 
    Then the barycentric coordinates admit the following representation in terms of oriented volumes:
    \[
    \lambda_i(x) = \frac{|K_i(x)|_{\mathrm{or}}}{|K|_{\mathrm{or}}}, \qquad i = 1, \dots, d+1,
    \]
    where $|\cdot|_{\mathrm{or}}$ denotes the oriented volume of the $d$ dimension. 
    As we specialize to one and two dimensions, these formulas degenerate to the familiar ratio formulas for lengths and areas, respectively. 
        
    For any point $x \in K$, its barycentric coordinates $\lambda_1(x), \dots, \lambda_{d+1}(x)$ satisfy
    \[
    \sum_{i=1}^{d+1} \lambda_i(x) = 1, \qquad  x = \sum_{i=1}^{d+1} \lambda_i(x)x^{a_i}.
    \]         
    Denote the vector $\bm{\lambda}_K=(\lambda_1,\dots,\lambda_{d+1})^\top$ and define the augmented vertex matrix by
    \[
    M_K=
    \begin{pmatrix}
    x^{a_1} & \cdots & x^{a_{d+1}}\\
    1 & \cdots & 1
    \end{pmatrix}.
    \]
    Let \(c_1,\dots,c_{d+1}\in\mathbb R^{d+1}\) be the columns of the adjugate matrix \(\operatorname{adj}(M_K)\), whose components are
    \[
    (c_j)_i
    =
    (-1)^{i+j}\det M_K^{(j,i)},
    \qquad i,j=1,\dots,d+1,
    \]
    where \(M_K^{(j,i)}\) is obtained from \(M_K\) by deleting its \(j\)th row and \(i\)th column. 
    Then the barycentric-coordinate map admits the unique affine representation
   \begin{equation}\label{eq:coord}
    \bm{\lambda}_K(x)=W_Kx+b_K,
    \end{equation} 
    where \(W_K\in\mathbb R^{(d+1)\times d}\) and
    \(b_K\in\mathbb R^{d+1}\) are given by
    \[
    W_K
    =
    \frac{1}{\det M_K}
    \begin{pmatrix}
    c_1 & \cdots & c_d
    \end{pmatrix},
    \qquad
    b_K  =
    \frac{1}{\det M_K}c_{d+1}.
    \]
    The definition \eqref{eq:coord} shows that, on each element, the barycentric coordinates are produced by a single-layer neural network with linear activation, whose weight matrix and bias vector are \(W_K\) and \(b_K\):
    \[
    x\mapsto \bm{\lambda}_K(x)=W_Kx+b_K.
    \]
Taken together, the local‑to‑global map \(g_K\) and the affine barycentric map \(\bm{\lambda}_K(\boldsymbol{x})=W_K\boldsymbol{x}+\boldsymbol{b}_K\) enable a direct encoding of finite element functions within the neural network.
The former characterizes the coupling of degrees of freedom across elements, whereas the latter embeds the geometric information of each simplex via network weights and biases. 
In what follows, we first construct the neural network formula for barycentric coordinates, then devise network architectures to represent linear finite element functions, and further generalize the proposed framework to a high ‐ order Lagrange finite element function. Lastly, we investigate the differentiability of the neural representations obtained. 

\subsection{Network architecture of Linear Lagrange Elements}
    Building on the local-to-global representation and the barycentric-coordinate formulation introduced above, we construct a mesh-induced sparsely connected network for linear Lagrange finite element functions. 
    For a fixed simplicial mesh, the network structure is determined directly by the finite element discretization, and the resulting network space reproduces the linear Lagrange finite element space exactly. 
    \begin{proposition}\label{pro3e1}[Exact neural representation of linear Lagrange finite elements]
    Let $\mathcal{T}_h=\{K_i\}_{i=1}^{c}$ be a simplicial mesh of    $\Omega\subset\mathbb{R}^d$, and let $V_h^1$ be the linear Lagrange finite element space associated with $\mathcal{T}_h$. 
    Then, for any $u_h\in V_h^1$, there is a mesh-induced sparsely connected neural network $u_\theta$ such that
    \[
    u_\theta(x)=u_h(x),
    \qquad
    \forall x\in\Omega.
    \]
    The corresponding mesh-level architecture is
    \[
    N_{\mathrm{nn}}^{(1)}
    =
    \left[
    d,\,
    (c(d+1),\sigma_1),\,
    (2c,\sigma_c),\,
    (4c,\sigma_2),\,
    1
    \right],
    \]
    where \((n,\sigma)\) denotes a hidden layer with \(n\) neurons and activation function \(\sigma\), respectively. 
    The explicit definitions of the activation functions $\sigma_1, \sigma_c, \sigma_2$ are given in the subsequent arguments. 
    \end{proposition}

    In the following, we verify Proposition \ref{pro3e1} by exactly constructing the explicit representation of a given linear finite element function $u_h\in V_h^1$. 

    For the linear simplicial Lagrange element considered in this section, the finite element space is
    \[
    V_h^1
    =
    \{v_h\in C^0(\Omega): v_h|_K\in P_1(K), \forall\, K\in\mathcal T_h\},
    \]
    where \(P_1(K)\) is the space of affine functions on \(K\). 
    The degrees of freedom are the nodal values at the vertices of the mesh, and hence $N_{\mathrm{dof}}=N$. 
    For each simplex \(K\in\mathcal T_h\), the linear Lagrange element has \(n_{\mathrm{loc}}=d+1\) local degrees of freedom given by the nodal values 
    \[
    \bm{u}_K=(u_{g_K(1)},\dots,u_{g_K(d+1)})^\top,
    \]
    and the corresponding shape functions coincide exactly with the barycentric coordinates of the element. 
    Then the finite element function is 
    \begin{equation}\label{eq:uk}
    \left.u_h\right|_K(x)
    =\sum_{i=1}^{d+1} u_{g_K(i)}\phi_i(x)=
    \bm{u}_K^\top(W_Kx+b_K).
    \end{equation}
    Hence, the local finite element evaluation has the form of a two-layer linear network: an element-geometry layer computes the basis values, and a nodal-value layer forms the final weighted superpose. 
    
\begin{figure}[htbp]
		\includegraphics[width=0.7\textwidth]{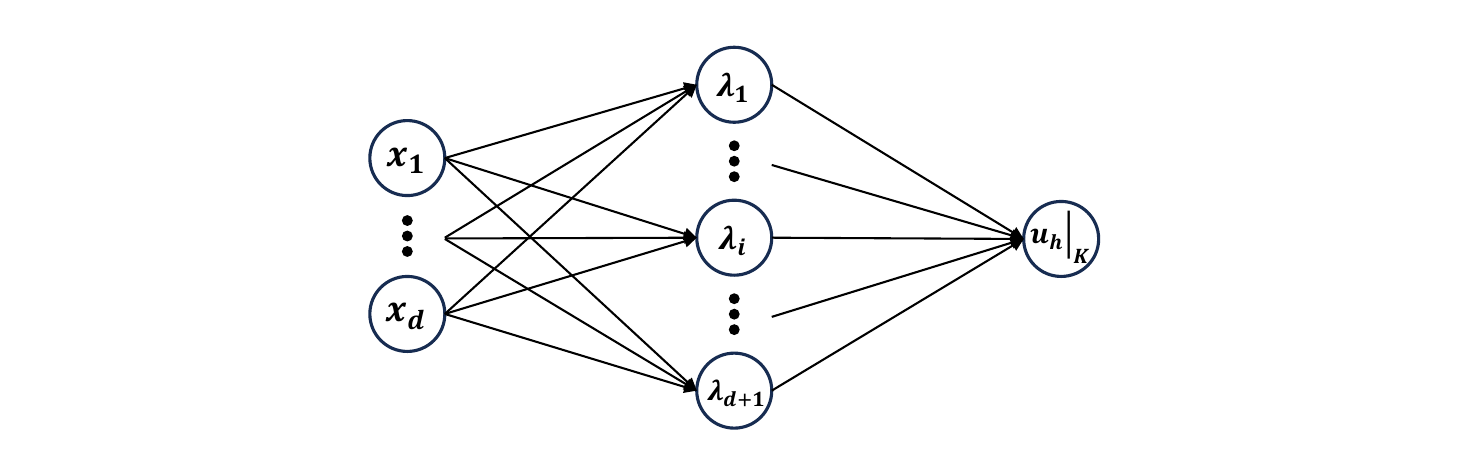}
		\caption{Local interpolation network with architecture  $N_{nn}=[d, d+1, 1]$ on a simplex.}
		\label{fig:localNET1}
	\end{figure}

    \Cref{fig:localNET1} depicts the local network structure for a scalar finite element function on a $d$-simplex $K$. 
    With architecture $N_{nn}=[d, d+1, 1]$, the input layer has $d$ neurons for the physical coordinate $x\in\mathbb R^d$, the hidden layer has $d+1$ neurons for the local basis values $\lambda_1(x),\dots,\lambda_{d+1}(x)$, and the output layer has one neuron for the elementwise value $\left.u_h\right|_K(x)$. 
    
    Equation \eqref{eq:uk} and \Cref{fig:localNET1} illustrate the neural-network representation of the functions of linear finite elements on element $K$. 
    This local constructions above give an exact realization of the finite element basis on each simplex. 
    Noting that functions represented by neural networks are defined globally over the domain $\Omega$, we construct an element identification network leveraging the intrinsic properties of barycentric coordinates to extend the local finite-element neural network representation onto the full computational domain. 
    Once the barycentric coordinates of \(x\) with respect to \(K\) have been calculated, their signs determine whether \(x\) lies in the interior of \(K\), on its boundary, or outside \(K\). 
    For a point \(x\in\mathbb R^d\), the sign pattern of \(\bm{\lambda}_K(x)\) determines its position relative to \(K\). 
    More precisely, the position of \(x\) relative to \(K\) is determined by the smallest barycentric coordinate:
    \[
    \min_{1\le i\le d+1}\lambda_i(x) 
    \ge 0, \quad {\textit{if}}\quad x\in K,\qquad {\textit{and}}\qquad 
    \min_{1\le i\le d+1}\lambda_i(x) <0, \quad {\textit{if}}\quad x\notin K. 
    \]
    Accordingly, each zero barycentric coordinate \(\lambda_i(x)\) identifies the boundary face opposite the vertex \(a_i\). 
    When several components vanish simultaneously, the point \(x\) lies where the corresponding boundary faces meet, which is a lower-dimensional face of \(K\). 
    Since the barycentric coordinates sum to one, at most $d$ components can vanish simultaneously, leaving the single remaining component equal to $1$. 
    However, in its direct form, it requires the check of all the barycentric-coordinate components \(d+1\) separately. 
    For a network realization, these componentwise checks are more naturally represented by a single scalar output indicating whether \(x\) belongs to the element. 
    Denote the ReLU activation by
    \[
    \sigma_1(t):=\operatorname{ReLU}(t)=\max\{t,0\}.
    \]
    Applying \(\sigma_1\) componentwise to the barycentric coordinates gives
    \[
    \lambda_i^*(x):=\sigma_1(\lambda_i(x)),
    \qquad i=1,\dots,d+1,
    \]
    which truncates the negative components to zero. 
    Consequently, the sum of the transformed barycentric coordinates remains equal to one for \(x\in K\) and becomes strictly larger than one for \(x\notin K\). 
    Thus the original componentwise identification test is reduced to a comparison between a single scalar quantity and the threshold value $1$, namely
    \[
    \sum_{i=1}^{d+1}\lambda_i^*(x)
    =1, \quad {\textit{if}}\quad x\in K,\qquad {\textit{and}}\qquad \sum_{i=1}^{d+1}\lambda_i^*(x)   >1, \quad {\textit{if}}\quad x\notin K.
    \]
    This comparison defines a nonlinear selection rule. 
    To encode it in the network, we introduce a custom activation function $\sigma_c(t)$, that is, \emph{cell activation},
    \[
    \sigma_c(t)=
    \begin{cases}
    1, & t\le 1,\\
    0, & t>1.
    \end{cases}
    \]
    In other words, \emph{cell activation} returns \(1\) on \(K\) precisely when \(x\in K\), and returns \(0\) otherwise. 

    The identification mechanism can be added as a parallel branch to the built local network structure \eqref{eq:uk}. 
    For this branch, the aggregation weights are given by the all-ones vector
    \[
    \mathbf 1_{d+1}:=(1,\dots,1)^\top\in\mathbb R^{d+1}.
    \]
    Once the ReLU-activated barycentric coordinates \(\bm{\lambda}_K^*(x)=(\lambda_1^*(x),\ldots,\lambda_{d+1}^*(x))^\top\) are available, the branch produces the element-selection coefficient
    \begin{equation}\label{eac}
    \alpha_K(x)
    =
    \sigma_c\!\left(\mathbf 1_{d+1}^\top\bm{\lambda}_K^*(x)\right).
    \end{equation}
    
\begin{figure}[htbp]
		\includegraphics[width=0.7\textwidth]{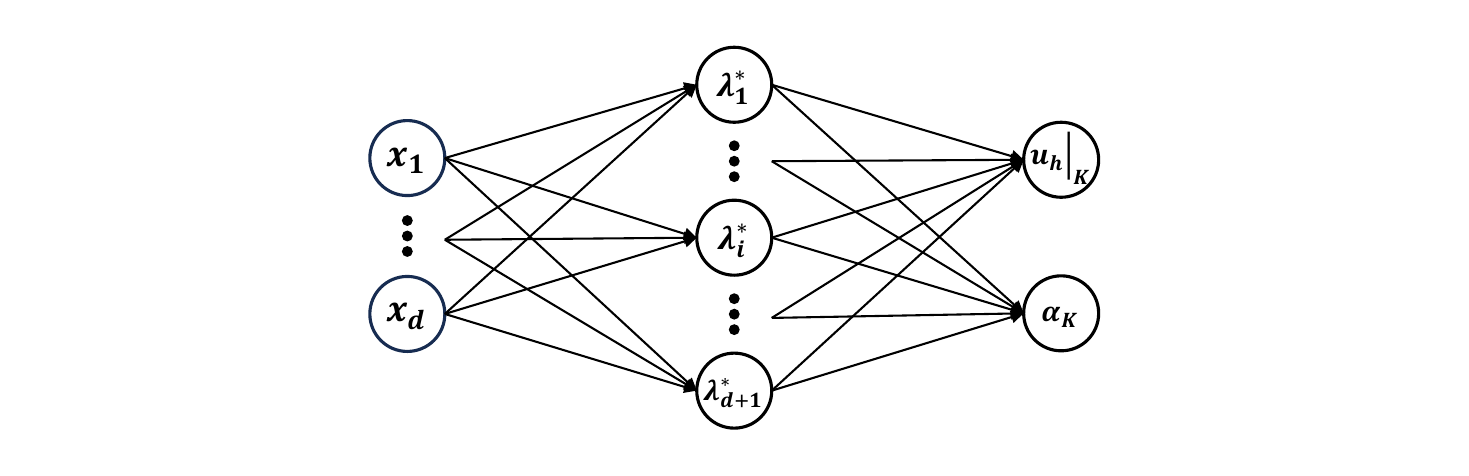}
		\caption{Local network with architecture $N_{nn}=[d, d+1, 2]$ for interpolation and element identification.}
		\label{fig:localNET2}
	\end{figure}

    The formula above shows that, on each element, the local element indicator coefficient is produced by a single-layer neural network with nonlinear activation. 
    Combined with the local finite element realization from the \eqref{eq:uk}, the identification mechanism leads to the two-output local network structure shown in \Cref{fig:localNET2}. 
    For a scalar finite element function on a simplex, the shared hidden layer of the local network with architecture $N_{nn}=[d,d+1,2]$ is given by
    \[
    \bm{\lambda}_K^*(x)
    =
    \operatorname{ReLU}(W_Kx+b_K).
    \]
    The hidden layer is shared by two output branches: one computes the local finite element value $\left.u_h\right|_K(x)$, and the other computes the element indicator coefficient $\alpha_K(x)$. 

    Having established neural representations for element-wise finite-element functions together with the element identification network, we proceed to discuss the construction of global neural-network approximations of finite-element functions over the entire computational domain. 
    For a query point \(x\), the blocks over \(\mathcal T_h\) produce two mesh-indexed vectors, namely the local finite element values \(\bigl(\left.u_h\right|_K(x)\bigr)_{K\in\mathcal T_h}\) and the element indicator coefficients \(\bigl(\alpha_K(x)\bigr)_{K\in\mathcal T_h}\). 
    The global sparsely connected neural network is obtained by multiplying these two mesh-indexed vectors componentwise and then summing over all elements:
    \begin{equation}\label{lenn}
    u_h(x)=\sum_{K\in\mathcal T_h}\alpha_K(x)\,\left.u_h\right|_K(x).
    \end{equation}
    It should be noted that suitable revisions must be applied to the expression of element-by-element of the finite element function \eqref{lenn} to accommodate scenarios with \(x\) located on the boundary of the element. 
    For this purpose, an averaging operator is introduced. Define the element patch of a point \(x\),
    \[
    \mathcal N_K(x):=\{K\in\mathcal T_h:\alpha_K(x)=1\},
    \]
    and set the contribution of indicator coefficients from element $K$ as 
    \[
    \widetilde\alpha_K(x):=A_K(\alpha_K(x))=\frac{1}{|\mathcal N_K(x)|}\alpha_K(x),
    \]
    where $A_K(\cdot)$ denotes the averaging operator and $|\mathcal N_K(x)|$ denotes the number of elements in $\mathcal{N}_K(x)$. 
    Accordingly, the assembled sparsely connected neural network output can be written in normalized form. 
    \begin{equation} 
        u_h(x)=\sum_{K\in\mathcal T_h}
        \widetilde{\alpha}_K(x)\,\left.u_h\right|_K(x).
        \label{eq:global-normalized-assembly}
    \end{equation}
        
\begin{figure}[htbp]
		\includegraphics[width=0.7\textwidth]{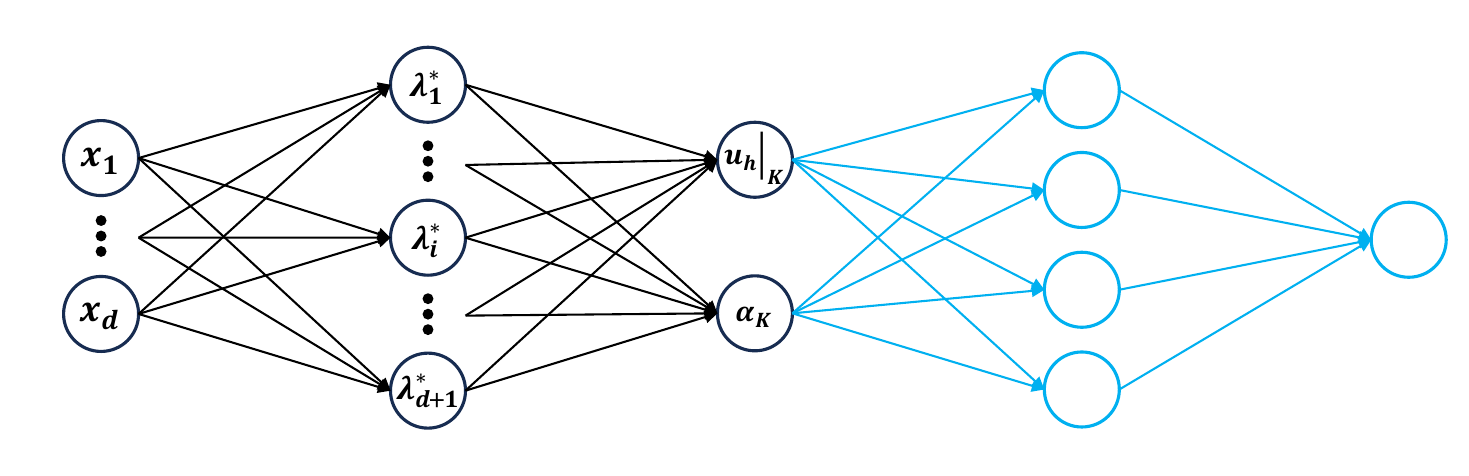}
		\caption{Local finite element network block with a product module. 
        }
		\label{fig:localNET3}
	\end{figure}

    However, the assembly formula \eqref{eq:global-normalized-assembly} requires the product of two same-layer outputs. 
    This multiplication is not a direct operation in a standard feedforward layer. 
    Following the construction in \cite{He2023}, we insert a product module that takes two scalar inputs \(a\) and \(b\) and returns their product \(ab\). 
    Define the squared ReLU activation
    \[
    \sigma_2(z):=\operatorname{ReLU}(z)^2 = \max\{z,0\}^2.
    \]
    Then, for any two scalars \(a,b\in\mathbb R\),
    \[\begin{aligned}
    ab
    =&
    \frac{1}{4}\bigl((a+b)^2-(a-b)^2\bigr)\\
    =&
    \frac{1}{4}\bigl(
    \sigma_2(a+b)+\sigma_2(-a-b)
    -\sigma_2(a-b)-\sigma_2(-a+b)
    \bigr),
    \end{aligned}\]
    where the identity \(\sigma_2(t)+\sigma_2(-t)=t^2\) has been used. 
    Thus componentwise multiplication can be realized exactly by one hidden layer using the squared ReLU activation:
    \begin{equation}\label{cw}
    ab
    =
    \begin{pmatrix}
    \mathbf 1_2^\top & -\mathbf 1_2^\top
    \end{pmatrix}
    \sigma_2\!\left(
    W
    \begin{pmatrix}
    a\\
    b
    \end{pmatrix}
    \right),
    \qquad
    W=
    \begin{pmatrix}
    \frac12 & \frac12\\
    -\frac12 & -\frac12\\
    \frac12 & -\frac12\\
    -\frac12 & \frac12
    \end{pmatrix}.
    \end{equation}
    Attaching the product module to the two-branch structure in \Cref{fig:localNET2} yields the product-coupled local network block $N_{nn}=[d, d+1, 2, 4, 1]$ shown in \Cref{fig:localNET3}. 
    For each element \(K\in\mathcal T_h\), its output can be written as
    \[
    \widetilde{\alpha}_K(x) \left.u_h\right|_K(x)
    =
    \begin{pmatrix}
    \mathbf 1_2^\top & -\mathbf 1_2^\top
    \end{pmatrix}
    \sigma_2\!\left(
    W
    \begin{pmatrix}
    u_K^\top \sigma_1(W_Kx+b_K)\\[2mm]
    A_K\left(\sigma_c\!\left(\mathbf 1_{d+1}^\top\sigma_1(W_Kx+b_K)\right)\right)
    \end{pmatrix}
    \right).
    \]

    The global network architecture is obtained by placing one such local block on each element of the mesh and coupling these blocks through the shared global degrees of freedom. 
    More precisely, the nodal coefficients associated with common mesh vertices are shared by all neighboring element blocks, while the element-selection coefficients determine which local contribution is active at a given input point. 
    The resulting architecture is therefore sparse at two levels: inside each element block, only the barycentric coordinates of that element are used; across the mesh, element blocks interact only through shared nodal coefficients and the final assembly. 

\begin{figure}[htbp]
    \centering
        \includegraphics[width=0.7\textwidth]{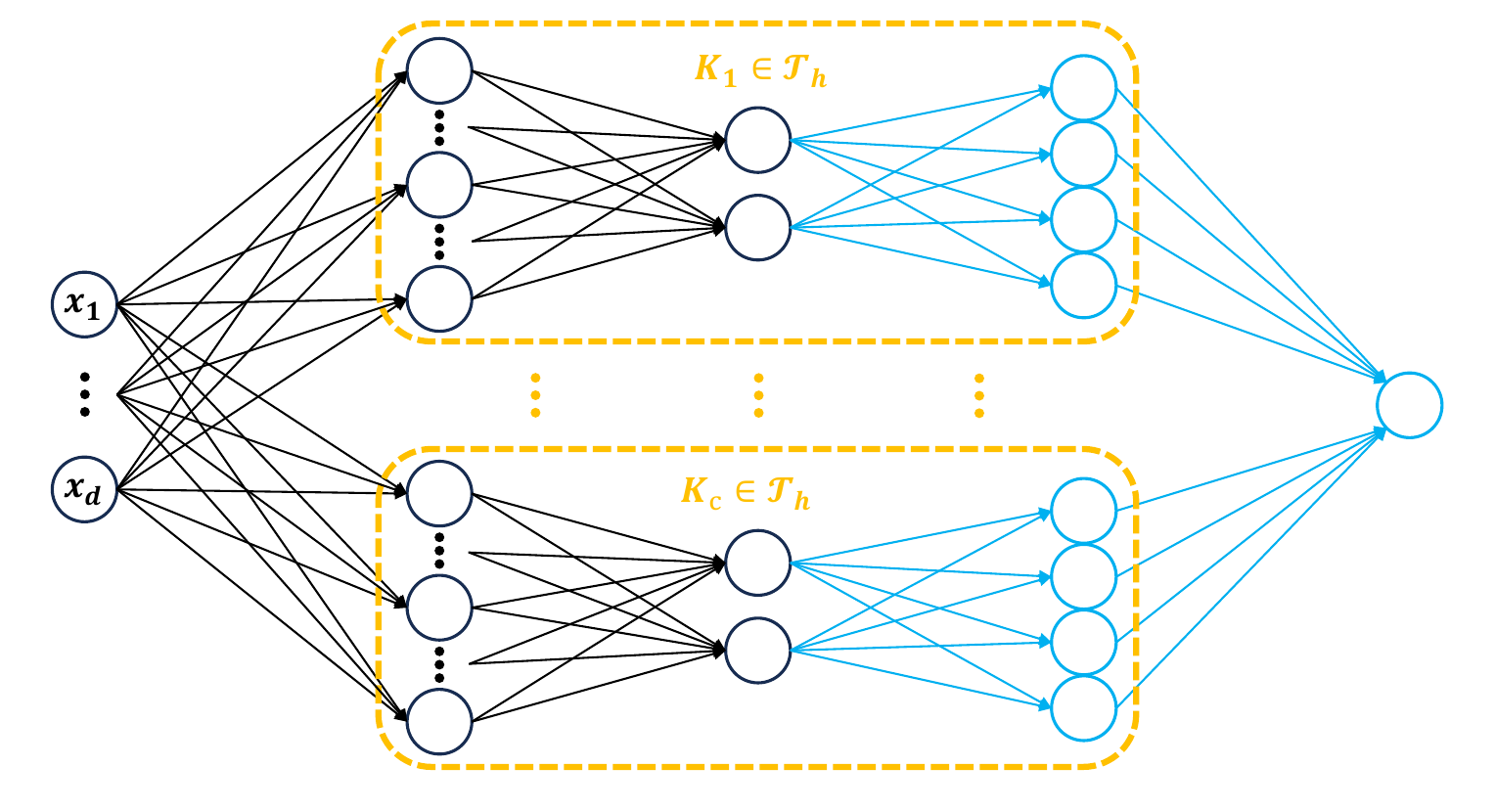}
        \caption{Sparsely connected neural network assembly for a linear finite element function.}
        \label{fig:Linear-network}
    \end{figure}

    \Cref{fig:Linear-network} shows the sparsely connected neural network of the linear finite element function in $V_h^1$ defined on a given mesh $\mathcal{T}_h$, where \(K_1,\ldots,K_c\in\mathcal T_h\) denote the indexed simplices in the mesh and \(c\) is the total number of elements. 
    Each simplex gives rise to a local affine barycentric layer and a local finite element readout, together with an element-selection branch. 
    These elementwise contributions are then multiplied by the corresponding selection coefficients and summed to produce the global finite element value. 
    This example show explicitly how the usual finite element assembly over all elements in $\mathcal{T}_h$ is represented as a sparse neural-network architecture. 
    
    \begin{example}\label{exam2e1} 
    In this illustrative example, an L-shaped computational domain is discretized by a triangulation \(\mathcal T_h\) consisting of 12 triangular elements, as plotted in \Cref{fig:Lshape-linear-network-A}. Based on this mesh \(\mathcal T_h\), we construct a dedicated neural network architecture for linear finite-element functions defined on \(\mathcal T_h\), as shown in \Cref{fig:Lshape-linear-network-B}. The network is composed of 12 element-level neural modules, where each module realizes the local finite-element representation on its corresponding element. 
    By superposing outputs from all local modules, the global sparsely connected neural network representation over the entire mesh is eventually attained.

    \begin{figure}[htbp]
    \centering
        \begin{subfigure}[b]{0.28\textwidth}
            \includegraphics[width=\linewidth]{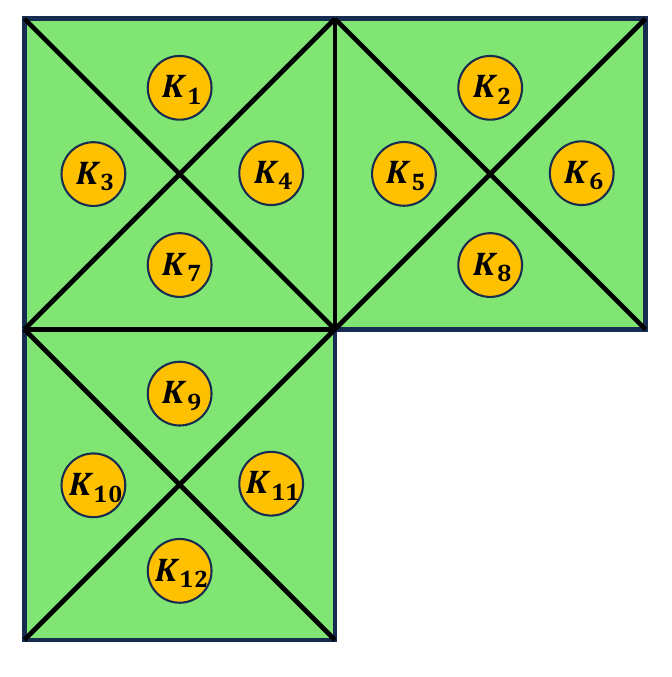}
            \caption{}
            \label{fig:Lshape-linear-network-A}
        \end{subfigure}
        \hspace{0.1cm}
        \begin{subfigure}[b]{0.7\textwidth}   
            \includegraphics[width=\linewidth]{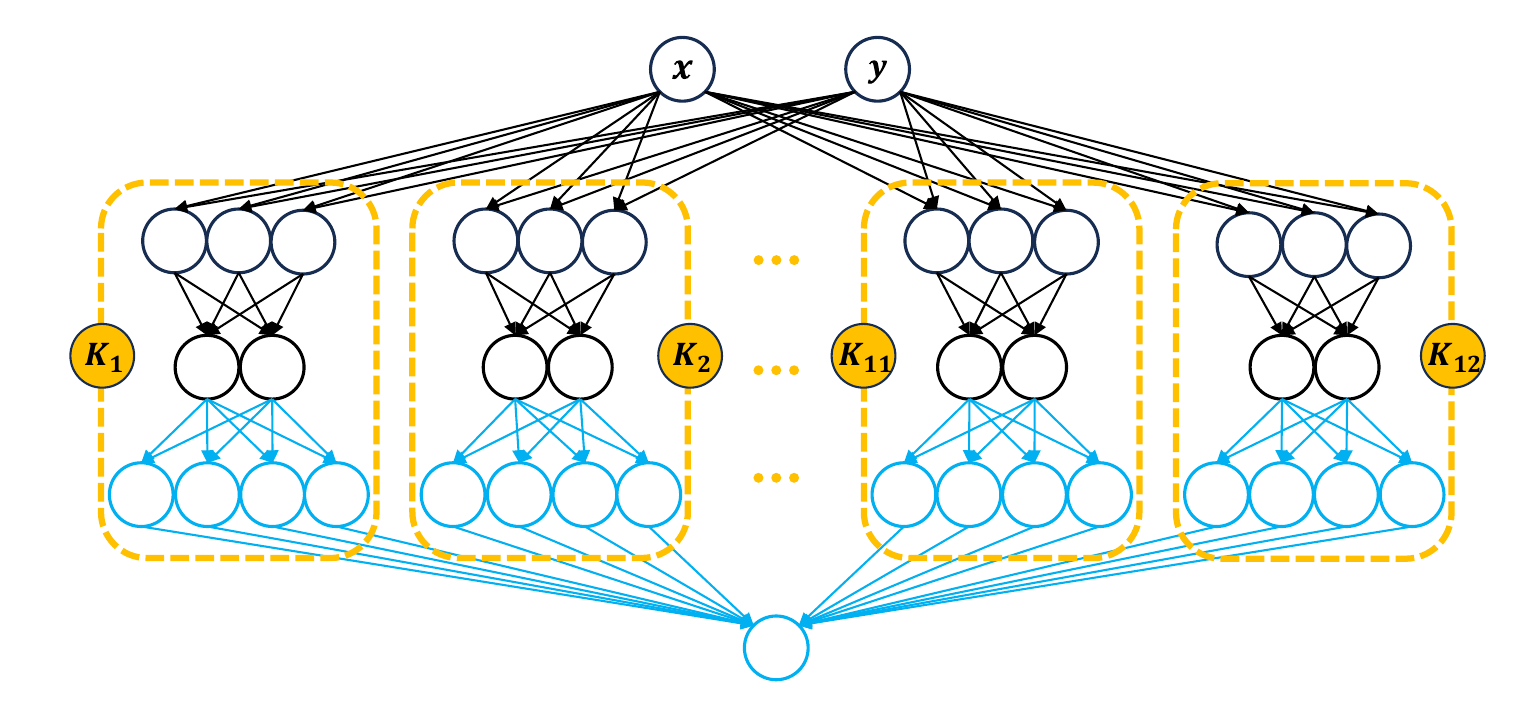}
            \caption{}
            \label{fig:Lshape-linear-network-B}
        \end{subfigure}
        \caption{Sparsely connected neural network assembly for a linear finite element function on an \(L\)-shaped mesh.}
        \label{fig:Lshape-linear-network}
    \end{figure}
    \end{example}

\subsection{Network architecture of Higher-Order Lagrange Elements}
This subsection extends the above construction to higher-order Lagrange element spaces
$V_h^p$ with $p\ge 2$.
\begin{proposition}\label{pro3e2}[Exact neural representation of $p$th-order Lagrange finite elements]
Let $\mathcal{T}_h=\{K_i\}_{i=1}^{c}$ be a simplicial mesh of
$\Omega\subset\mathbb{R}^d$, and let $V_h^p$ be the $p$th-order
Lagrange finite element space associated with $\mathcal{T}_h$.
Then, for any $u_h\in V_h^p$, there is a mesh-induced sparsely
connected neural network $u_\theta$ such that
\[
u_\theta(x)=u_h(x),
\qquad
\forall x\in\Omega .
\]
The corresponding mesh-level architecture is
\[
N_{\mathrm{nn}}^{(p)}
=
\left[
d,\,
\bigl(c(d+1),\sigma_1\bigr),\,
\bigl(c\,2^p n_{\mathrm{loc}},\sigma_p\bigr),\,
c\,n_{\mathrm{loc}},\,
(2c,\sigma_c),\,
\bigl(4c,\sigma_2\bigr),\,
1
\right],
\]
where \((n,\sigma)\) denotes a hidden layer with \(n\) neurons and the
activation function \(\sigma\), respectively. The activation function $\sigma_p$ is defined in \eqref{eqp}, and $\sigma_1, \sigma_2, \sigma_c$ are the same as in Proposition \ref{pro3e1}.
\end{proposition}

We now prove this proposition using a constructive argument. Following the same local-to-global route as before, we begin with the definition of the local $p$th-order Lagrange basis on a single simplex element. 
    
    Let $K$ be a nondegenerate $d$-simplex with affinely independent vertices $a_1,\dots,a_{d+1}$. 
    Accordingly, a $p$th-order Lagrange element on $K$ is determined by $n_{\mathrm{loc}}$ local interpolation points together with the nodal interpolation conditions.    
    Compared with the linear case, which uses only the \(d+1\) vertex interpolation points, a higher-order element contains many more local interpolation points, already \(\binom{3+2}{2}=10\) for a two-dimensional cubic element. 
    On simplices, a uniform parametrization is obtained by indexing interpolation points via barycentric multi-indices of order $p$, formally defined as
    \begin{equation}\label{eq2e2}
    \mathcal I_p^d
    :=
    \{\bm\nu=(\nu_1,\dots,\nu_{d+1})\in\mathbb N_0^{d+1}:|\nu|=p\}.
    \end{equation}
    Each \(\bm\nu\in\mathcal I_p^d\) determines an interpolation point \(a_{\bm\nu}\in K\), defined by
    \[
    a_{\bm\nu}=
    \sum_{m=1}^{d+1}\frac{\nu_m}{p} x^{a_m}.
    \]
    Since \(\frac{|\bm\nu|}{p}=1\), the above affine representation identifies the barycentric coordinate vector of \(a_{\bm\nu}\) as
    \[
    \bm{\lambda}_K(a_{\bm\nu})
    =
    \left(
    \frac{\nu_1}{p},\dots,\frac{\nu_{d+1}}{p}
    \right)^\top .
    \]
    With the nodal points written in barycentric form, the construction of the corresponding basis functions can be expressed directly in terms of the barycentric coordinates. 
    For each $\bm\nu\in\mathcal I_p^d$, the local basis function $\phi_{\bm\nu}\in P_p(K)$ is characterized by the condition
    \[
    \phi_{\bm\nu}(a_{\bm\mu})=\delta_{\bm\nu \bm\mu},
    \qquad \bm\mu\in\mathcal I_p^d.
    \]
    Equivalently, $\phi_{\bm\nu}$ takes the value $1$ at its associated interpolation point $a_{\bm\nu}$ and vanishes at all other interpolation points $a_{\bm\mu}$ with $\bm\mu\neq\bm\nu$.            
    The following standard formula gives these basis functions explicitly in barycentric coordinates. 
    \begin{lemma}\label[lemma]{lem:HlocalNET}\cite{Nicolaides1972}
        Let $K$ be the simplex, $\bm{\lambda}_K$ its barycentric-coordinate vector, 
        $\mathcal I_p^d$ the multi-index set, and 
        $\{a_{\bm\nu}\}_{\bm\nu\in\mathcal I_p^d}$ the corresponding interpolation points defined above. 
        For each $\bm\nu\in\mathcal I_p^d$, the local $p$th-order Lagrange basis function
        $\phi_{\bm\nu}\in P_p(K)$ associated with $a_{\bm\nu}$ satisfies
        \[
        \phi_{\bm\nu}(x)
        =
        \prod_{m=1}^{d+1}
        \prod_{r=0}^{\nu_m-1}
        \frac{p\,\lambda_m(x)-r}{\nu_m-r},
        \qquad x\in K.
        \]
        Here the product over \(r\) is interpreted as \(1\) when \(\nu_m=0\).
    \end{lemma}

    \Cref{lem:HlocalNET} shows that each higher-order local basis function is obtained by multiplying affine factors of the barycentric coordinates. 
    Let $\mathcal I_p^d=\{\bm\nu_1,\dots,\bm\nu_{n_{\mathrm{loc}}}\}$, each local basis function can be represented as a polynomial function of the barycentric-coordinate vector:
    \[
    \phi_{\bm\nu_i}(x)
    =
    q_i(\bm{\lambda}_K(x)),
    \qquad i=1,\dots,n_{\mathrm{loc}},
    \]
    where \(q_i:\mathbb R^{d+1}\to\mathbb R\) is determined by the product formula in \Cref{lem:HlocalNET}. 
    Collecting these functions gives the local basis vector
    \[
    \bm{\phi}_K(x)
    =
    \bigl(
    q_1(\bm{\lambda}_K(x)),
    \dots,
    q_{n_{\mathrm{loc}}}(\bm{\lambda}_K(x))
    \bigr)^\top .
    \]
    Thus the local higher-order basis evaluation follows the chain 
    \[
    x
    \longmapsto
    \bm{\lambda}_K(x)
    \longmapsto
    \bm{\phi}_K(x).
    \]
    The first layer is the affine barycentric layer
    \[
    \bm{\lambda}_K(x)=W_Kx+b_K,
    \]
    which is exactly the same geometric layer as in the linear element. 
    The higher-order part then forms degree-\(p\) combinations of the \(d+1\) barycentric coordinates. 
    Following the product-module construction provided above, each \(q_i\) is implemented by a \(p\)-fold product module, whose internal activation is
    \begin{equation}\label{eqp}
    \sigma_p(z):=\operatorname{ReLU}(z)^p = \max\{z,0\}^p.
    \end{equation}
    Consequently, all local \(p\)th-order basis values are produced in parallel from \(\bm{\lambda}_K(x)\), giving a basis layer with \(n_{\mathrm{loc}}\) outputs. 
    The linear case is recovered when \(p=1\), as this layer reduces to the affine barycentric layer itself and each \(q_i\) simply selects one barycentric-coordinate component. 
    \begin{figure}[htbp]
    \centering
        \includegraphics[width=0.5\textwidth]{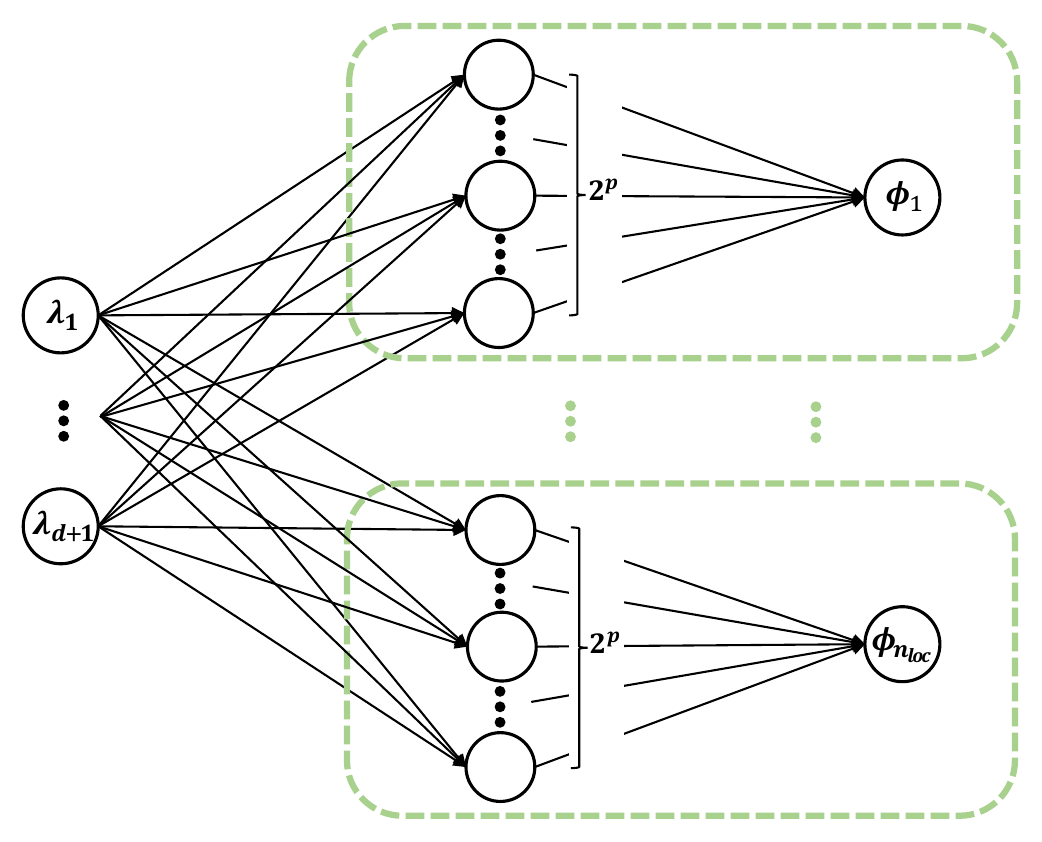}
        \caption{Local basis-realization network for a \(p\)th-order Lagrange element. 
        }
        \label{fig:HlocalNET1}
    \end{figure}
    
    \Cref{fig:HlocalNET1} presents the basis-realization block for the map 
    \(\bm{\lambda}_K(x)\mapsto\bm{\phi}_K(x)\). 
    Starting from the \(d+1\) barycentric coordinates, which coincide with the linear local basis values, the hidden layer uses \(p\)-fold product modules with internal activation \(\sigma_p(z)\) to generate higher-order basis functions. 
        The intermediate representation has \(2^p n_{\mathrm{loc}}\) neurons and feeds into \(n_{\mathrm{loc}}\) outputs, one for each local \(p\)th-order basis value.
    This block is also a sparsely connected subnetwork whose connection pattern is prescribed by the barycentric factors appearing in the polynomials \(q_i\). 
    The hidden layer is arranged into \(n_{\mathrm{loc}}\) groups of \(2^p\) neurons, where group \(i\) implements the \(\sigma_p\)-based product representation of \(q_i\) and contributes only to the output component \(\phi_{\bm\nu_i}\). 
    With \(\bm{\phi}_K(x)\) available as the output of the basis-realization layer, the elementwise finite element value is formed by the standard nodal readout:
    \[
    \left.u_h\right|_K(x)
    =
    \bm{u}_K^\top\bm{\phi}_K(x)
    =
    \sum_{\bm\nu\in\mathcal I_p^d}
    u_{g_K(\bm\nu)}\phi_{\bm\nu}(x),
    \qquad x\in K,
    \]
    where $g_K(\bm\nu)$ denotes the same local-to-global map in multi-index form, namely $g_K(\bm\nu_i):=g_K(i), i=1,\dots,n_{\mathrm{loc}}$.

In parallel to the extraction of the nodal output, each local module also produces a scalar-valued element-selection coefficient \(\alpha_K(x)\).
Consistent with the linear finite-element neural network defined in  \eqref{lenn}, the element-selection coefficient \(\alpha_K(x)\) is calculated by activating the inner-product against the summed high-order basis functions. 
Denote the high-order basis vector evaluated at the ReLU-activated barycentric coordinates
\(\bm{\lambda}_K^*(x)\) by
    \[
    \bm{\phi}_K^*(x)
    =
    \bigl(
    q_1(\bm{\lambda}_K^*(x)),
    \dots,
    q_{n_{\mathrm{loc}}}(\bm{\lambda}_K^*(x))
    \bigr)^\top.
    \]
    The element-selection coefficient is then given by
    \[
    \alpha_K(x)
    =
    \sigma_c\!\left(
    \mathbf 1_{n_{\mathrm{loc}}}^{\top}\bm{\phi}_K^*(x)
    \right).
    \]
    Using the generalized Vandermonde identity for the barycentric product basis, we obtain
    \[
    \mathbf 1_{n_{\mathrm{loc}}}^{\top}\bm{\phi}_K^*(x)
    =
    \sum\limits_{\nu\in\mathcal I_p^d}
    \prod_{m=1}^{d+1}
    \binom{p\,\lambda_{K,m}^*(x)}{\nu_m}
    =
    \binom{
    p\sum\limits_{m=1}^{d+1}\lambda_{K,m}^*(x)
    }{p}
    =
    \prod_{r=0}^{p-1}
    \frac{
    p\sum\limits_{m=1}^{d+1}\lambda_{K,m}^*(x)-r
    }{
    p-r
    }.
    \]
    Since \(\sum\limits_{m=1}^{d+1}\lambda_{K,m}^*(x)=1\) for \(x\in K\) and is larger than \(1\) for \(x\notin K\), the right-hand side is equal to \(1\) in the former case and larger than \(1\) in the latter case.
    The higher-order basis-realization output admits the matrix form
    \begin{equation}\label{cwb}
    \bm{\phi}_K^*(x)
    =
    R_p
    \sigma_p\!\left(
    W_{d,p}\sigma_1(W_Kx+b_K)+b_{d,p}
    \right),
    \end{equation}
    where \(W_{d,p}\) and \(b_{d,p}\) encode the affine factors in the \(p\)th-order barycentric product representations, and \(R_p\) is the block-sparse readout matrix
    \begin{equation}\label{cr}
    R_p
    =
    I_{n_{\mathrm{loc}}}
    \otimes r_p
    \in \mathbb R^{n_{\mathrm{loc}}\times 2^p n_{\mathrm{loc}}}.
    \end{equation}
    Here \(r_p\in\mathbb R^{1\times 2^p}\) is the readout row of a single \(p\)-fold product module. 
    Combining the local finite element output and the element-selection coefficient by the same product module as before gives the local product-coupled higher-order block
    \begin{equation}\label{cww}
    \widetilde{\alpha}_K(x)\,\left.u_h\right|_K(x)
    =
    \begin{pmatrix}
    \mathbf 1_2^\top & -\mathbf 1_2^\top
    \end{pmatrix}
    \sigma_2\!\left(
    W
    \begin{pmatrix}
    \bm{u}_K^\top\bm{\phi}_K^*(x)\\[2mm]
    A_K\left(\sigma_c\!\left(
    \mathbf 1_{n_{\mathrm{loc}}}^\top\bm{\phi}_K^*(x)
    \right)\right)
    \end{pmatrix}
    \right),
    \end{equation}
where $W$ is defined in \eqref{cw}, and the global sparsely connected neural network for the finite element function is 
\begin{equation}\label{cnfe}
u_h(x)=\sum\limits_{K\in\mathcal{T}_h}\widetilde{\alpha}_K(x)\,\left.u_h\right|_K(x).
\end{equation}    
    This completes the elementwise construction for \(p\)th-order Lagrange elements. 
    The subsequent mesh-level assembly follows the same pattern as in the linear finite element neural network \eqref{eq:global-normalized-assembly}: the elementwise outputs are selected by the normalized element indicators and then summed over all elements. 
    Neighboring elements are coupled through the common global coefficients associated with shared Lagrange degrees of freedom. 
    Hence the higher-order case introduces only a new local basis-realization block, while the global coupling and assembly structure remain unchanged. 

Subsequently, we construct neural network representations for quadratic finite element on interval and cubic finite element over triangular element, respectively. Only the local neural network formulations on individual element is explicitly provided herein, the corresponding global neural representation can be recovered via superposition following \eqref{cnfe}.

\begin{example}\label{exam2e2}
    We consider the quadratic element in one dimensional. 
    For \(d=1\) and \(p=2\), the index set \(\mathcal I_2^1\) consists of the three multi-indices
    \[
    (2,0),\qquad (0,2),\qquad (1,1).
    \]
    The corresponding local quadratic basis functions on an interval are
    \[
    \phi_{(2,0)}=\lambda_1(2\lambda_1-1),
    \qquad
    \phi_{(0,2)}=\lambda_2(2\lambda_2-1),
    \qquad
    \phi_{(1,1)}=4\lambda_1\lambda_2.
    \]
    For the basis functions associated with the vertex, they share the same one-variable representation:
    \[
    \lambda_i(2\lambda_i-1)
    =
    \sigma_2\!\left(\frac{3\lambda_i-1}{2}\right)
    +
    \sigma_2\!\left(\frac{1-3\lambda_i}{2}\right)
    -
    \sigma_2\!\left(\frac{1-\lambda_i}{2}\right)
    -
    \sigma_2\!\left(\frac{\lambda_i-1}{2}\right),
    \qquad i=1,2.
    \]
    For the basis function associated with the midpoint, the product identity gives
    \[
    4\lambda_1\lambda_2
    =
    \sigma_2(\lambda_1+\lambda_2)
    +
    \sigma_2(-\lambda_1-\lambda_2)
    -
    \sigma_2(\lambda_1-\lambda_2)
    -
    \sigma_2(-\lambda_1+\lambda_2).
    \]
    \Cref{fig:exam1D} shows the resulting local network block for this quadratic element. 
    After the barycentric-coordinate layer, the three identities above are organized as three parallel squared-ReLU modules, highlighted in blue in the figure, with \(n_{\mathrm{loc}}=3\) and \(p=2\). 
    The resulting basis values are then used to form two intermediate scalar quantities, the elementwise finite element value \(\left.u_h\right|_K(x)\) and the element-selection coefficient \(\alpha_K(x)\), which are subsequently coupled by a product module. 

    \begin{figure}[htbp]
    \centering
        \includegraphics[width=0.8\textwidth]{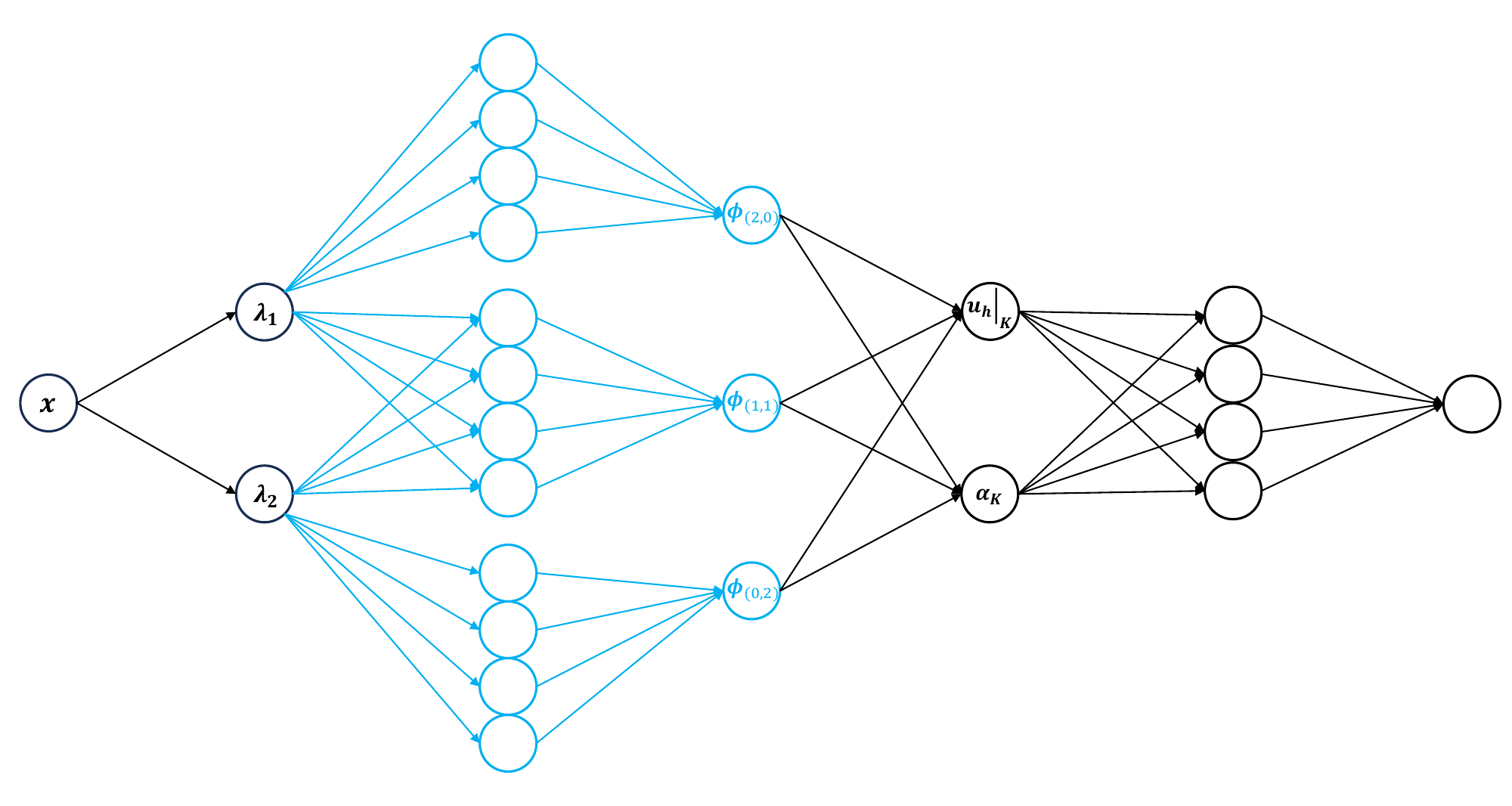}
        \caption{Local network realization for the $1$d quadratic finite element function on element $K$.
        }
        \label{fig:exam1D}
    \end{figure}    
\end{example}

\begin{example}\label{exam2e3}
    For a two-dimensional triangular element
    \[
    K=\mathrm{conv}\{a_1, a_2, a_3\}\subset\mathbb R^2,
    \]
    with
    \[
    a_1=(x_1,y_1)^\top,\qquad a_2=(x_2,y_2)^\top,\qquad a_3=(x_3,y_3)^\top,
    \] 
    the barycentric coordinates of $P=(x,y)^\top\in K$ are
    \[
    \lambda_1(P)=\frac{S_{\triangle PBC}}{S_{\triangle ABC}},\qquad
    \lambda_2(P)=\frac{S_{\triangle PCA}}{S_{\triangle ABC}},\qquad
    \lambda_3(P)=\frac{S_{\triangle PAB}}{S_{\triangle ABC}},
    \]
    where
    \[
    S_{\triangle ABC}
    :=
    (x_2-x_1)(y_3-y_1)-(y_2-y_1)(x_3-x_1)
    \]
    denotes the oriented area factor of the triangle $ABC$.

    Hence
    \[
    \begin{bmatrix}
    \lambda_1\\
    \lambda_2\\
    \lambda_3
    \end{bmatrix}
    =
    W_K
    \begin{bmatrix}
    x\\
    y
    \end{bmatrix}
    +b_K,
    \]
    with
    \[
    W_K=
    \frac{1}{S_{\triangle ABC}}
    \begin{bmatrix}
        y_2-y_3 & x_3-x_2\\
        y_3-y_1 & x_1-x_3\\
        y_1-y_2 & x_2-x_1
    \end{bmatrix},
    \qquad
    b_K=
    \frac{1}{S_{\triangle ABC}}
    \begin{bmatrix}
        x_2y_3-x_3y_2\\
        x_3y_1-x_1y_3\\
        x_1y_2-x_2y_1
    \end{bmatrix}.
    \]
 
    For the cubic Lagrange element on a triangle, as \(d=2\) and \(p=3\) in \eqref{eq2e2}. 
    The index set \(\mathcal I_3^2\) consists of the $10$ multi-indices
    \[
    (3,0,0),\quad (0,3,0),\quad (0,0,3),\quad
    (2,1,0),\quad (1,2,0),\quad (0,2,1),\quad (0,1,2),\quad
    (1,0,2),\quad (2,0,1),\quad
    (1,1,1).
    \]
    Using this ordering, the corresponding local basis functions are
    \begin{align*}
    \phi_1 &= \frac12(3\lambda_1-1)(3\lambda_1-2)\lambda_1, &
    \phi_2 &= \frac12(3\lambda_2-1)(3\lambda_2-2)\lambda_2, &
    \phi_3 &= \frac12(3\lambda_3-1)(3\lambda_3-2)\lambda_3, \\
    \phi_4 &= \frac92\,\lambda_1\lambda_2(3\lambda_1-1), &
    \phi_5 &= \frac92\,\lambda_1\lambda_2(3\lambda_2-1), &
    \phi_6 &= \frac92\,\lambda_2\lambda_3(3\lambda_2-1), \\
    \phi_7 &= \frac92\,\lambda_2\lambda_3(3\lambda_3-1), &
    \phi_8 &= \frac92\,\lambda_3\lambda_1(3\lambda_3-1), &
    \phi_9 &= \frac92\,\lambda_3\lambda_1(3\lambda_1-1), \\
    \phi_{10} &= 27\,\lambda_1\lambda_2\lambda_3.
    \end{align*}
    The cubic products in these formulas are realized by the three-variable product module. 
    With
    \[
    \sigma_3(t):=\operatorname{ReLU}(t)^3,
    \]
    one has
    \[
    \begin{aligned}
    abc
    = \frac1{24}\bigl[&
    \sigma_3(a+b+c)
    -\sigma_3(-a-b-c)
    +\sigma_3(-a+b-c)
    -\sigma_3(a-b+c) \\
    +&\sigma_3(a-b-c)
    -\sigma_3(-a+b+c)
    +\sigma_3(-a-b+c)
    -\sigma_3(a+b-c)
    \bigr].
    \end{aligned}
    \]
    For example, the basis function \(\phi_1\), associated with the vertex $a_1$, is obtained by substituting
    \[
    a=3\lambda_1-1,\qquad
    b=3\lambda_1-2,\qquad
    c=\lambda_1,
    \]
    which gives
    \[
    \begin{aligned}
    \phi_1
    =
    \frac1{48}\bigl[&
    \sigma_3(7\lambda_1-3)
    -\sigma_3(3-7\lambda_1)
    +\sigma_3(-\lambda_1-1)
    -\sigma_3(\lambda_1+1) \\
    +&\sigma_3(1-\lambda_1)
    -\sigma_3(\lambda_1-1)
    +\sigma_3(3-5\lambda_1)
    -\sigma_3(5\lambda_1-3)
    \bigr].
    \end{aligned}
    \]
    Similarly, for the basis function \(\phi_4\) associated with a edge, by substituting
    \[
    a=\lambda_1,\qquad
    b=\lambda_2,\qquad
    c=3\lambda_1-1
    \]
    yields
    \[
    \begin{aligned}
    \phi_4
    =
    \frac{3}{16}\bigl[&
    \sigma_3(4\lambda_1+\lambda_2-1)
    -\sigma_3(1-4\lambda_1-\lambda_2)
    +\sigma_3(1-4\lambda_1+\lambda_2) 
    -\sigma_3(4\lambda_1-\lambda_2-1) \\
    +&\sigma_3(1-2\lambda_1-\lambda_2)
    -\sigma_3(2\lambda_1+\lambda_2-1) 
    +\sigma_3(2\lambda_1-\lambda_2-1)
    -\sigma_3(1-2\lambda_1+\lambda_2)
    \bigr].
    \end{aligned}
    \]
    For the basis function \(\phi_{10}\), which is associated with the element center, by substituting
    \[
    a=\lambda_1,\qquad
    b=\lambda_2,\qquad
    c=\lambda_3
    \]
    gives
    \[
    \begin{aligned}
    \phi_{10}
    =
    \frac{9}{8}\bigl[&
    \sigma_3(\lambda_1+\lambda_2+\lambda_3)
    -\sigma_3(-\lambda_1-\lambda_2-\lambda_3)
    +\sigma_3(-\lambda_1+\lambda_2-\lambda_3)
    -\sigma_3(\lambda_1-\lambda_2+\lambda_3) \\
    +&\sigma_3(\lambda_1-\lambda_2-\lambda_3)
    -\sigma_3(-\lambda_1+\lambda_2+\lambda_3)
    +\sigma_3(-\lambda_1-\lambda_2+\lambda_3)
    -\sigma_3(\lambda_1+\lambda_2-\lambda_3)
    \bigr].
    \end{aligned}
    \]
    The remaining cubic basis functions are obtained by the same direct substitution of their three affine barycentric factors. 
    
    \Cref{fig:local2Dcubic} plots the neural network architecture of the cubic basis functions on a triangular element. 
    The selected branches for \(\phi_1\), \(\phi_4\), and \(\phi_{10}\) illustrate the basis function associated with vertex, edge, and interior, respectively. 
    In detail, starting from the input coordinates \((x,y)\), the network first computes the barycentric coordinates \(\lambda_1,\lambda_2,\lambda_3\).
    The followed branch is a \(\sigma_3\)-based basis-realization module: the hidden neurons evaluate the affine barycentric combinations in the cubic representation, and the output node forms the corresponding fixed linear combination. 
    The figure shows three representative local cubic basis functions, namely the vertex basis function \(\phi_1\), the edge basis function \(\phi_4\), and the interior basis function \(\phi_{10}\). The remaining cubic local basis functions are obtained by the same construction.
    \begin{figure}[htbp]
    \centering
        \includegraphics[width=0.8\textwidth]{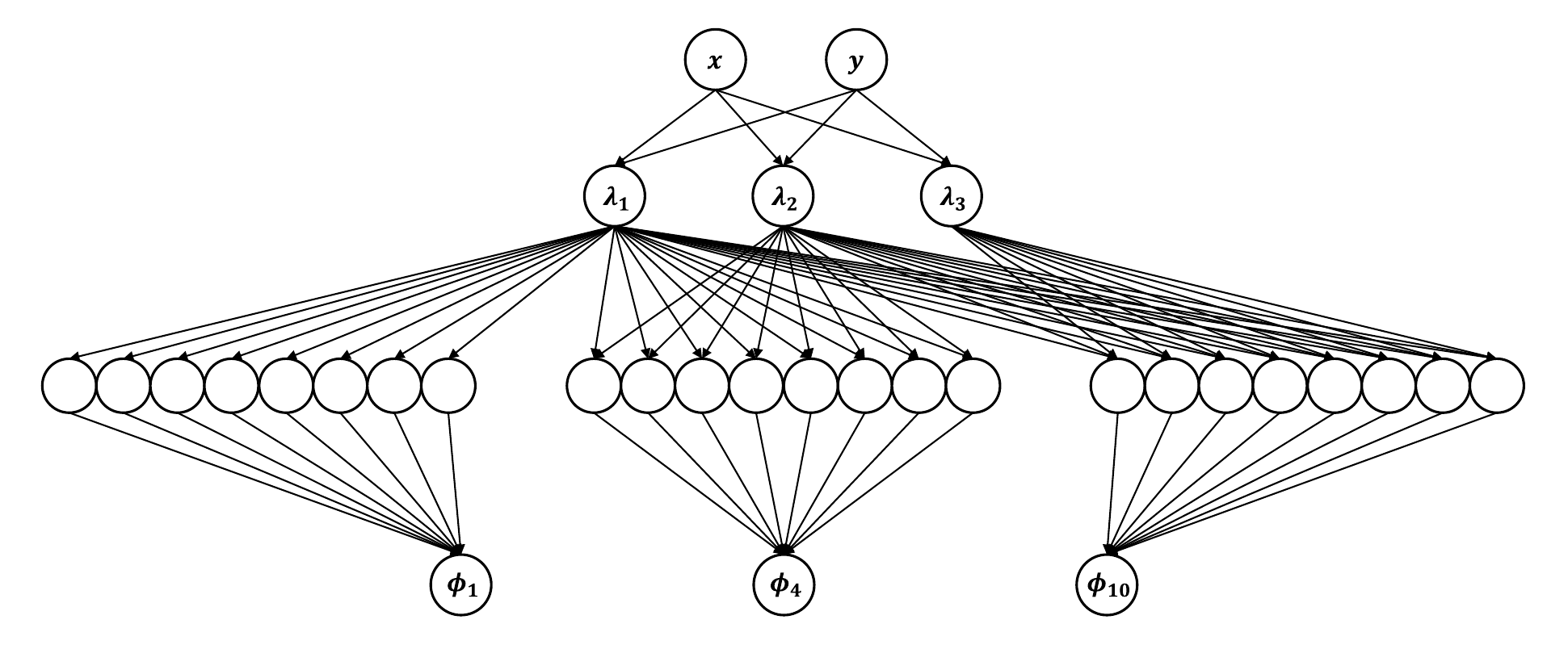}
        \caption{Neural network architecture of the cubic basis functions on a triangular element.}
        \label{fig:local2Dcubic}
    \end{figure}  
     
\end{example}


\subsection{Differentiability}
    In the foregoing subsections, we establish a sparsely connected neural network representation for Lagrange finite element functions, which integrates barycentric affine mappings, local basis function implementation modules, element selection coefficients, product operation modules, and global assembly operations. The constructed neural network space is strictly equivalent to the Lagrange finite element space, such that all inherent analytical properties of finite element spaces, including classical error estimation theories, can be fully inherited by the proposed network framework. Beyond inherent approximation capabilities, the developed architecture supports explicit derivative computation and treats mesh vertices as trainable optimization variables. This distinctive feature enables the construction of a fully differentiable mesh adaptation framework for finite element analysis.

    In a standard finite element implementation, both the value and the derivative of \(u_h\) at a point \(x\) are evaluated from the same local basis functions and their gradients:
    \[
    \left.u_h\right|_K(x)=\sum_i u_i\,\phi_i(x),\qquad
    \nabla_x \left.u_h\right|_K(x)=\sum_i u_i\,\nabla_x\phi_i(x).
    \]
    While the two quantities rely on the same set of coefficients \(u_i\), they are typically computed via separate loops or stored separately.
    By contrast, the sparsely connected neural network representation unifies these two evaluations within a single computational graph: the forward pass directly yields \(u_h(x)\), and differentiation of the same graph yields \(\nabla_x u_h(x)\). 
    
    Applying the chain rule to \eqref{eq:global-normalized-assembly} gives the global gradient
    \begin{equation}\label{eq:global-gradient}
        \nabla_x u_h(x)
        =
        \sum_{K\in\mathcal T_h}
        \left(
        \nabla_x\widetilde{\alpha}_K(x)\,\left.u_h\right|_K(x)
        +
        \widetilde{\alpha}_K(x)\,\nabla_x \left.u_h\right|_K(x)
        \right).
    \end{equation}
    The first term requires the gradient of the normalized indicator. 
    By differentiating the normalization relation, $\nabla_x\widetilde{\alpha}_K$ is expressed entirely in terms of the backward derivatives of the unnormalized indicators $\alpha_{K'}$:
    \[
    \nabla_x\widetilde{\alpha}_K(x)
    =
    \frac{
    \left(\sum\limits_{K'}\alpha_{K'}(x)\right)
    \nabla_x\alpha_K(x)
    -
    \alpha_K(x)
    \sum\limits_{K'}\nabla_x\alpha_{K'}(x)
    }{
    \left(\sum\limits_{K'}\alpha_{K'}(x)\right)^2
    }.
    \]
    Formal differentiation of the cell activation gives
    \[
    \nabla_x\alpha_K(x)
    =
    \sigma_c'(s_K(x))\,\nabla_x s_K(x).
    \]
    where \(s_K(x)\) is the scalar input associated with element \(K\).
    The cell activation lies on constant branches whenever $s_K(x) \neq 1$, yielding $\sigma_c'(s_K(x)) = 0$, whereas the classical derivative remains undefined at the threshold $s_K(x) = 1$. 
    Similar to ReLU-type activations, this nondifferentiable point requires a suitable backward derivative assignment. 
    Since both adjacent branches of \(\sigma_c\) possess zero derivative, we impose the consistent threshold value:
    \[
    \sigma_c'(1):=0.
    \]
    With this backward convention, the normalization of the element-selection coefficients gives
    \[
    \nabla_x\widetilde{\alpha}_K(x)=0.
    \]
    
    As for $\nabla_x \left.u_h\right|_K(x)$, it is governed by the ReLU-activated barycentric map $\bm{\lambda}_K^*(x)$. 
    Analogously, the backward derivative of $\sigma_1$ presents an issue precisely at points where $W_Kx+b_K$ vanishes, i.e., on the faces of $K$. 
    Standard automatic-differentiation backends, such as PyTorch, conventionally assign $\sigma_1'(0)=0$, which treats a zero barycentric component as inactive in the backward pass. 
    For the elementwise barycentric map used here, however, a zero entry of \(W_Kx+b_K\) corresponds to a boundary value of an active element, and its backward derivative is taken from the active side. 
    We therefore modify the backward derivative to
    \[
    \sigma_1^*(t)=\sigma_1(t)=\max\{t,0\},
    \qquad
    (\sigma_1^*)'(t)
    =
    \begin{cases}
    0, & t<0,\\
    1, & t\ge 0,
    \end{cases}
    \]
    which changes the backward value at $t=0$ from $0$ to $1$. 
    Consequently, the derivative of the barycentric map becomes
    \[
    \nabla_x \bm{\lambda}_K^*(x) = (\sigma_1^*)'(W_Kx+b_K) W_K,
    \]
    retaining the full gradient contribution across element boundaries. 
    The zero backward contribution of the normalized element-selection coefficients gives
    \[
    \nabla_x u_h(x)
    =
    \sum_{K\in\mathcal T_h}
    \widetilde{\alpha}_K(x)\,
    \nabla_x \left.u_h\right|_K(x),
    \]
    where the modified activation \(\sigma_p^*\) is used, satisfying
    \[
    \sigma_p^*(t)=\operatorname{ReLU}(t)^p,\qquad (\sigma_p^*)^{(k)}(t)=\frac{p!}{(p-k)!}\operatorname{ReLU}(t)^{p-k}\,(\sigma_1^*)'(t),\quad 1\le k\le p.
    \]
With the backward-mode derivatives predefined for \(\sigma_c\) and \(\sigma_p^*\), the gradient \(\nabla_x u_h(\bm{x})\) can be readily evaluated via automatic differentiation applied directly on the established sparse computational graph. This obviates the need for standalone finite-element gradient subroutines, as the derivatives of local basis functions are naturally yielded through network backpropagation.
    

\section{Application}\label{secApp}
In this section, we first present a numerical example to illustrate that the constructed neural network architecture is capable of accurately representing Lagrange finite elements. 
Furthermore, taking advantage of the mesh-free feature of this neural network representation, finite element functions can be interpolated onto new meshes. This strategy is incorporated into adaptive finite element algorithms for parabolic equations \cite{Hao2026}.

\begin{example}[Accuracy test for the neural network representation of finite element functions] 
    We test the accuracy of the proposed neural network representation for $d\in\{1,2,3\}$. 
    Let $\mathcal{T}_h=\{K_i\}_{i=1}^{N_K}$ denote the Delaunay mesh generated from a prescribed set of uniformly distributed points, with $\{z_i\}_{i=1}^{N_V}$ and $\{c_i\}_{i=1}^{N_E}$ denoting the mesh vertices and edge midpoints, respectively, and $\{g_j^{K_i}\}_{j=1}^{m}$ denoting the $m$ Gauss quadrature points on each element $K_i$, $i=1,\ldots,N_K$. $V_h$ represents the linear finite element space over $\mathcal{T}_h$. For a given function 
    \[u=\left\{\begin{aligned}
        &\sin(x), \qquad & d=1,\\
        &\sin(x)\cos(y), \qquad & d=2,\\
        &\sin(x)\cos(y)e^{z}, \qquad & d=3,
    \end{aligned}\right.\]
    let $I_hu\in V_h$ denotes its finite element interpolation, and $u_\theta(x)$ is the neural network surrogate of $I_hu$. We calculate the errors of $e=u_\theta(x)-I_hu(x)$ as follows:
    \[\begin{aligned}
        E_{max}=&\max\limits_{s_i}\{|u_\theta(s_i)-I_hu(s_i)|\},\\
        E_{mean}=&\frac{\sum\limits_{s_i}|u_\theta(s_i)-I_hu(s_i)|}{N_s},
    \end{aligned}\]
    where $s_i=z_i$, $c_i$, or $g_j^{K_i}$, and $N_s$ is the number of total mesh vertices, edge midpoints or Gauss quadrature points, respectively.
    
    \Cref{tab:error_statistics_dimensions} summarizes the resulting error statistics. 
    For $d=1$, the mesh is constructed with $5000$ uniformly spaced points; for $d=2$ and $d=3$, Delaunay triangulation and tetrahedralization are performed on uniform point grids of size $50\times50$ and $10\times10\times10$, respectively. 
    For all values of $d$ and all sampling locations, the maximum errors are of order $10^{-15}$--$10^{-14}$, while the mean errors are generally of order $10^{-16}$ or smaller. 
    Comparable error levels are obtained at the nodes, edge midpoints, and Gaussian quadrature points. 
    In particular, the results at the non-nodal points show that direct evaluation through the neural network reproduces the values obtained by the conventional procedure of element location followed by local finite element interpolation.
    These results are consistent with the exact representation property established above, with the remaining discrepancies attributable to floating-point round-off.
\end{example}

\begin{table}[htbp]
\centering
\caption{Errors evaluated at different sampling points for $d=1,2,3$.}
\label{tab:error_statistics_dimensions}
\sisetup{
    scientific-notation = true,
    table-format = 1.1e-1,
    exponent-product = \times,
    retain-zero-exponent = true
}
\renewcommand{\arraystretch}{1.15}
\setlength{\tabcolsep}{10pt}

\begin{tabular}{clSSS}
    \toprule
    \textbf{Dimension}
    & \textbf{Error}
    & {\textbf{Nodes}}
    & {\textbf{Edge midpoints}}
    & {\textbf{Gaussian points}} \\
    \midrule

    \multirow{2}{*}{1D}
    & $E_{max}$ & 5.814793e-15 & 1.162959e-14 & 2.314815e-14 \\
    & $E_{mean}$ & 1.549726e-18 & 6.876712e-18 & 2.714969e-17 \\
    \addlinespace

    \multirow{2}{*}{2D}
    & $E_{max}$ & 5.218048e-15 & 1.010303e-14 & 9.992007e-15 \\
    & $E_{mean}$ & 3.757457e-16 & 4.339106e-16 & 7.363282e-16 \\
    \addlinespace

    \multirow{2}{*}{3D}
    & $E_{max}$ & 2.664535e-15 & 4.440892e-15 & 3.108624e-15 \\
    & $E_{mean}$ & 1.941727e-16 & 2.243627e-16 & 2.543556e-16 \\
    \bottomrule
\end{tabular}
\end{table}

\begin{figure}[htbp]
\centering
    \includegraphics[width=0.16\textwidth]{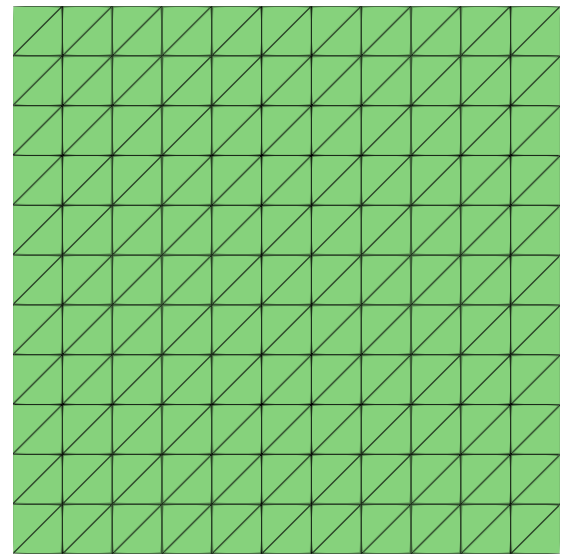}
    \includegraphics[width=0.16\textwidth]{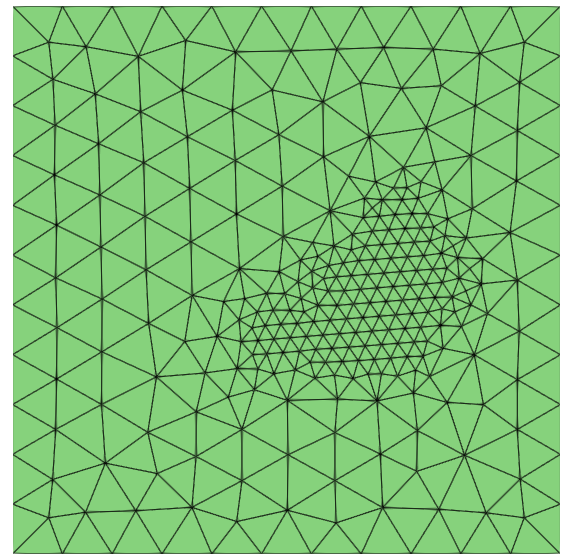}
    \includegraphics[width=0.16\textwidth]{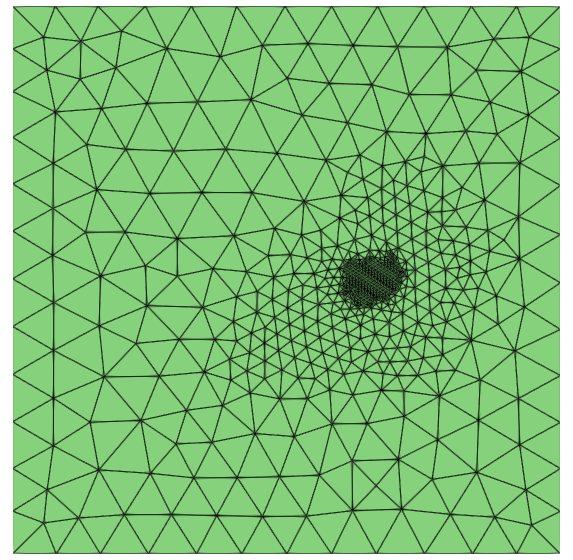}
    \includegraphics[width=0.16\textwidth]{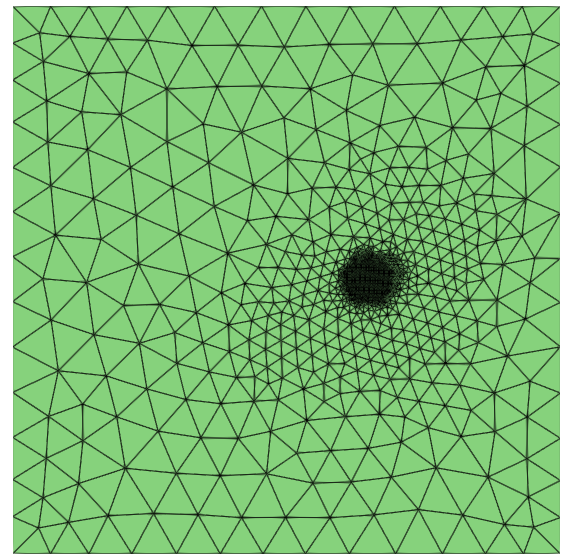}
    \includegraphics[width=0.16\textwidth]{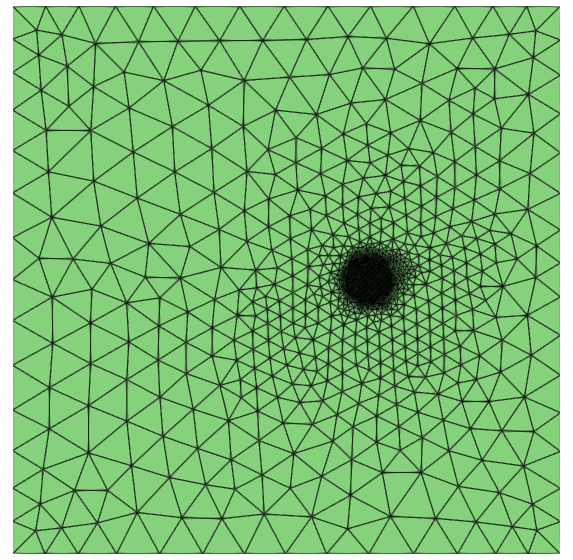}
    \includegraphics[width=0.16\textwidth]{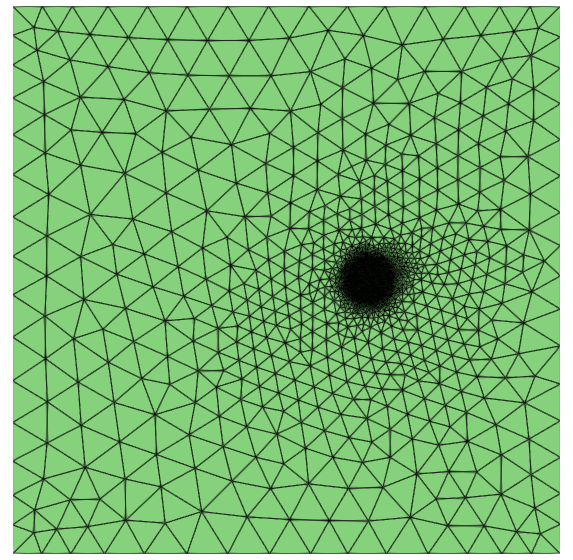}
    \includegraphics[width=0.16\textwidth]{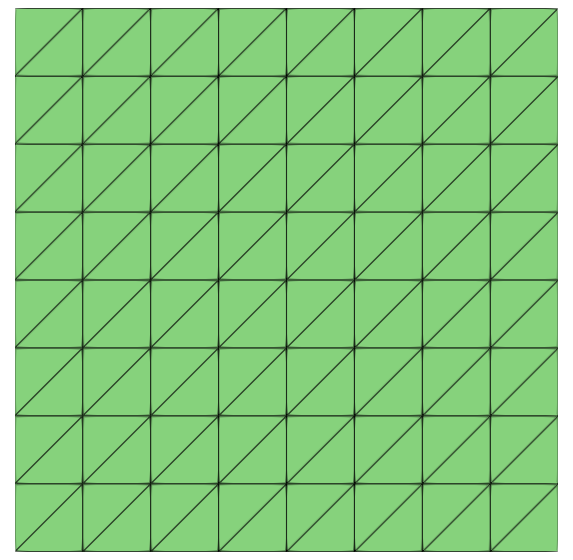}
    \includegraphics[width=0.16\textwidth]{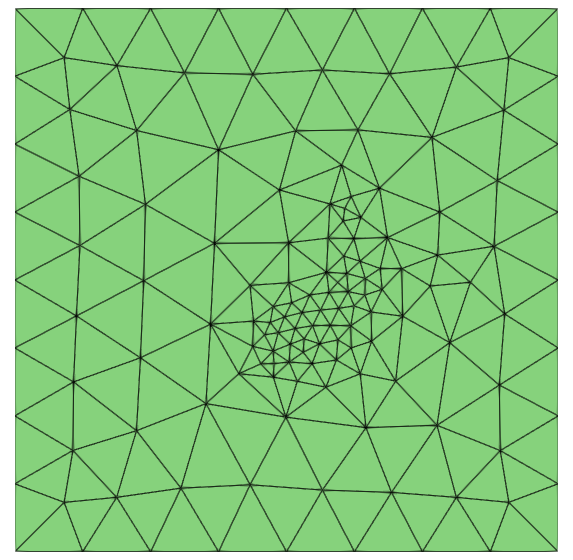}
    \includegraphics[width=0.16\textwidth]{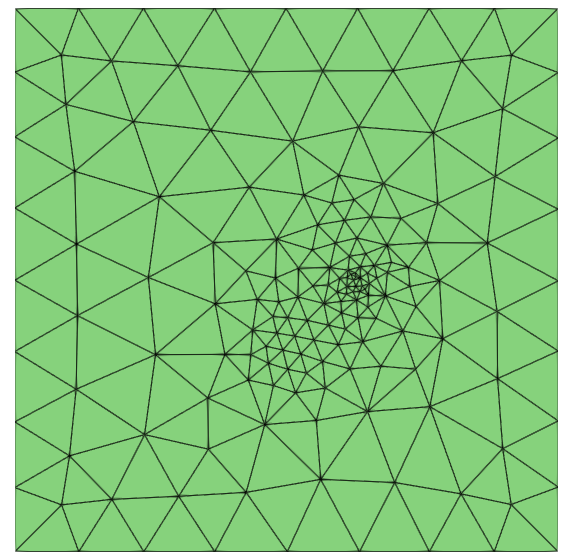}
    \includegraphics[width=0.16\textwidth]{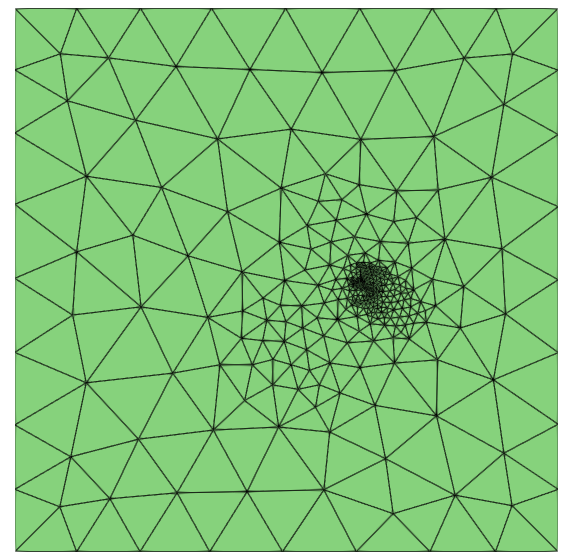}
    \includegraphics[width=0.16\textwidth]{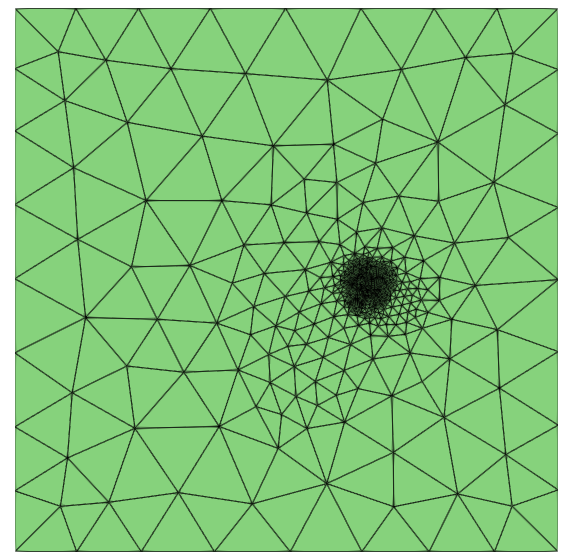}
    \includegraphics[width=0.16\textwidth]{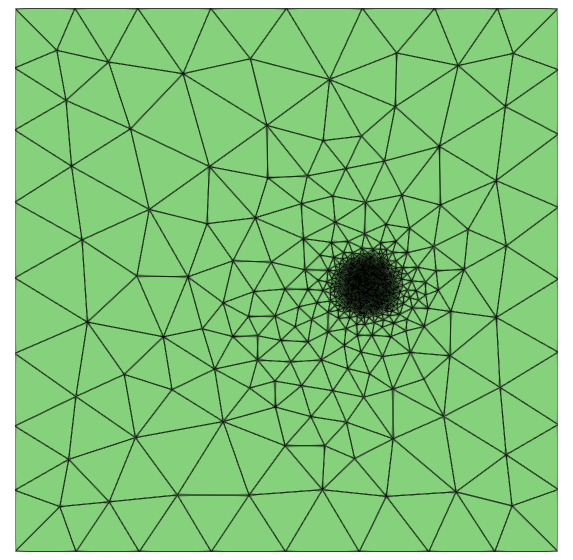}
    \caption{Successive adaptive meshes at the time $t=1.0$. Top: $P^1$-AFEM, bottom: $P^2$-AFEM.}
    \label{fig:HRFEM}
\end{figure}

\begin{example}[Neural network representation of FEM function in adaptive FEM] 
    We consider the neural network-enhanced adaptive finite element method for parabolic equations in~\cite{Hao2026}. 
    Let $\Omega=[-1,1]^2$ and consider 
    \begin{equation*}
    \left\{\begin{aligned} 
        u_t-\Delta u &= f, && \text{in } \Omega\times(0,1],\\ 
        u &= g, && \text{on } \partial\Omega\times(0,1],\\ 
        u(\bm{x},0) &= u_0, && \text{in } \Omega, 
    \end{aligned}\right.
    \end{equation*}
    where $f$, $g$, and $u_0$ are chosen such that the exact solution is 
    $$ u(x,y,t) = \exp\!\left[-500\bigl(x-0.3\cos(2\pi t)\bigr)^2\right] \exp\!\left[-500\bigl(y-0.3\sin(2\pi t)\bigr)^2\right]. $$
    In the original algorithm, a general neural network is trained to represent the finite element solution during the adaptive computation. 
    Here, this neural network is replaced by the exact finite element neural representation developed in this work, while the other components of the adaptive algorithm are kept unchanged. 
\end{example}
\begin{figure}[htbp]
\centering
    \includegraphics[width=0.35\textwidth]{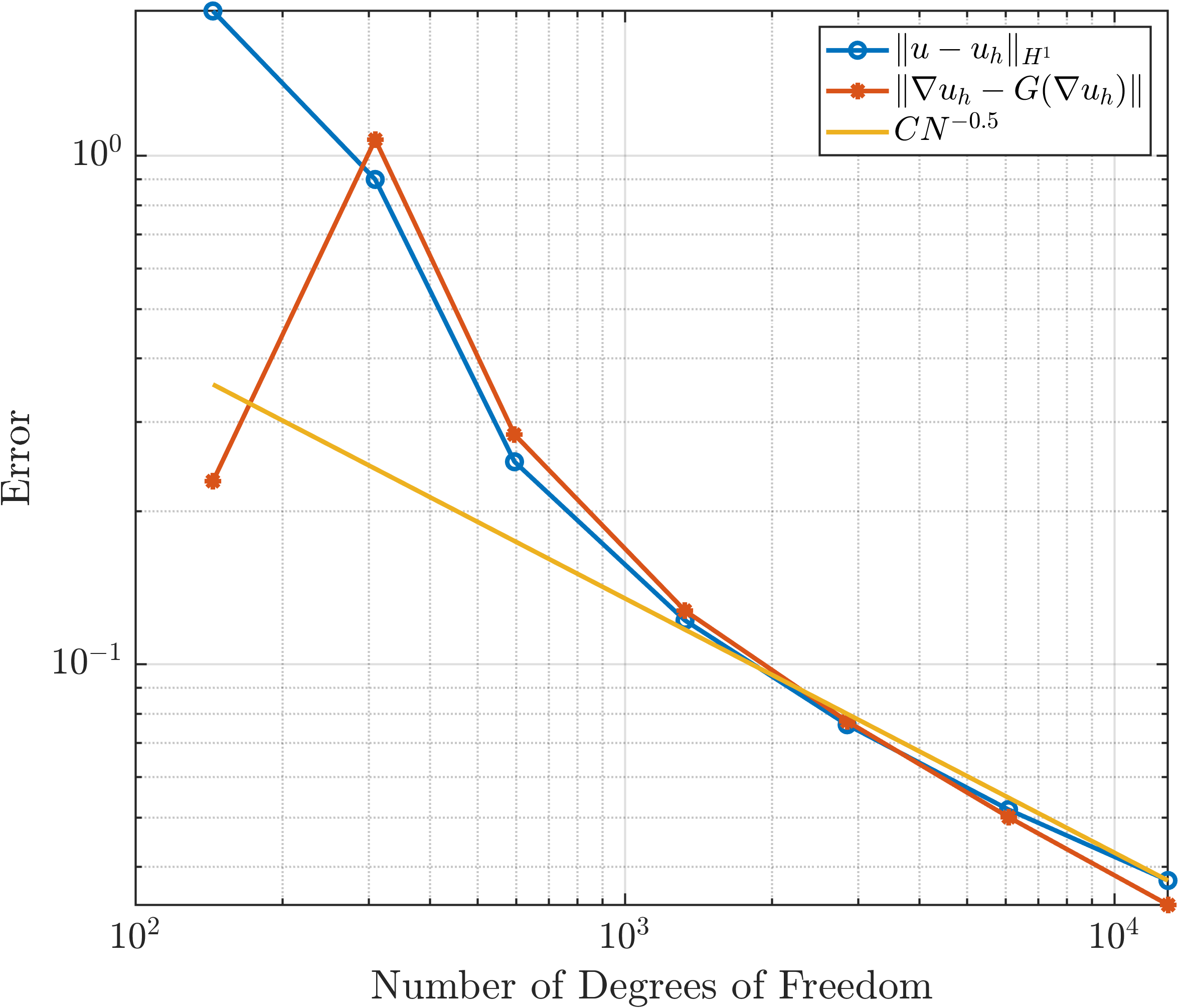}\hspace{2cm}
    \includegraphics[width=0.35\textwidth]{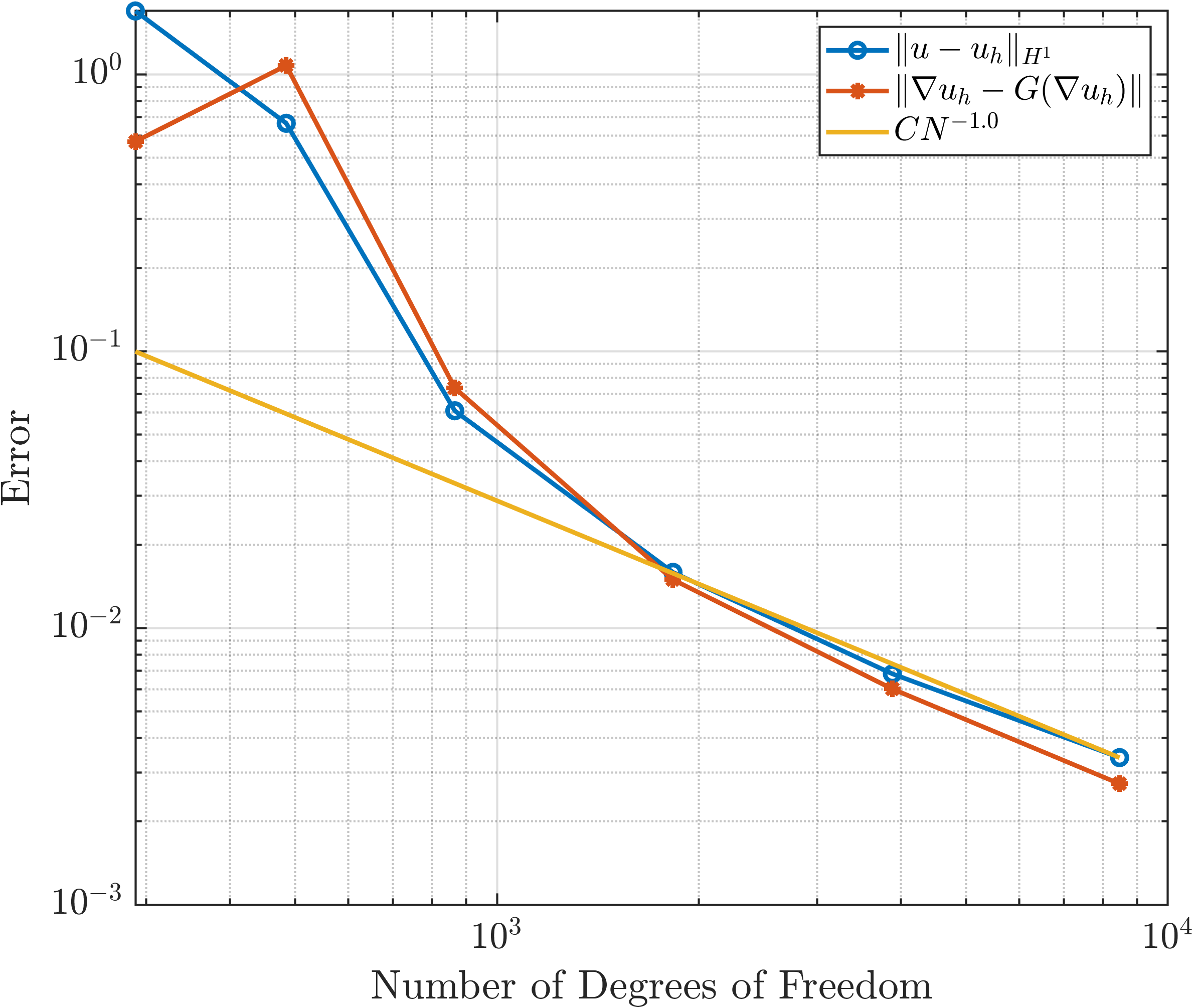}
    \caption{Convergence of the $H^1$ error and the recovery-based error estimator for the adaptive FEM using the neural representation proposed in this work: $P_1$ (left) and $P_2$ (right).}
    \label{fig:NFEM}
\end{figure}
\Cref{fig:NFEM} presents the convergence results for the $P_1$ and $P_2$ discretizations. 
The $H^1$ errors decay with the expected rates $N^{-1/2}$ and $N^{-1}$, respectively, where $N$ denotes the number of degrees of freedom. 
These rates agree well with the standard approximation orders of the corresponding finite element spaces in two dimensions. 
The corresponding adaptive mesh evolution over six successive adaptive steps at $t=1.0$ is shown in \Cref{fig:HRFEM}. 
For both finite element spaces, the refinement remains concentrated in the region where additional resolution is required, with little refinement pollution in the surrounding area. 
Some mesh vertices are also relocated during the adaptive procedure, so the meshes at successive steps are generally non-nested. 
Therefore, in the adaptive computation of finite elements in time $t^{n+1}$, the finite element solution obtained at $t^n$ needs to be interpolated onto the newly adapted mesh in the current time step. To ensure high efficiency and interpolation accuracy, this work adopts the neural network representation of finite element functions developed to achieve accurate solution representation.
Compared with the original approach, in which a general neural network must be trained to represent the finite element solution during the adaptive computation, the proposed representation is obtained directly from the finite element degrees of freedom and requires no training. 
Thus, the expected finite element convergence rates are recovered without the neural-network training step required in the original adaptive procedure.

\section{Conclusion}\label{secCon}
This work develops a complete framework for realizing finite element spaces via mesh-induced sparsely connected neural networks. For linear Lagrange elements, local basis functions are exactly equivalent to element barycentric coordinates, and the corresponding local network modules are constructed through explicit affine mappings. For high-order Lagrange elements and general polynomial basis functions, local finite element bases can be uniformly formulated as polynomial combinations of barycentric coordinates. Combined with degree-of-freedom matching and global assembly operations, the resulting network function space strictly matches the standard finite element space.

The proposed construction provides a novel perspective to interpret finite element spaces beyond classical numerical analysis, redefining them as a family of network function spaces with explicit geometric origins, well-defined parameter semantics, and inherent sparse connectivity. Since the established network space is theoretically identical to the conventional finite element space, classical finite element results, including interpolation error bounds, convergence rates, function enrichment properties, and mesh-dependent theoretical conclusions, can be directly inherited by the presented network framework.


In future investigations, we will extend the proposed neural representation framework to various other finite element spaces including edge finite element spaces. Using the computational merits inherent to neural networks, we will apply the established framework to $r$-adaptive mesh optimization, where mesh nodes are optimized as trainable parameters.

\section*{Acknowledgments}
	This research was supported by the National Key R$\&$D Program of China (2024YFA1012600) and NSFC Project (12431014).
	
\bibliographystyle{plain}
\bibliography{references}
	
\end{document}